\documentclass[11pt]{amsart}
\usepackage{amsmath, amsthm, amsfonts, amssymb} 
\usepackage{url}

\newtheorem{theo}{Theorem}[section]
\newtheorem{cor}[theo]{Corollary}

\newtheorem{prop}[theo]{Proposition}
\newtheorem{claim}[theo]{Claim}
\theoremstyle{definition}
\newtheorem{rem}[theo]{Remark}
\newtheorem{exam}[theo]{Example}

\newtheorem{df}[theo]{Definition}

\newcommand{\R}{\mathbf{R}}
\newcommand{\C}{\mathbf{C}}
\newcommand{\Z}{\mathbf{Z}}
\newcommand{\F}{\mathbf{F}}

\newcommand{\N}{\mathbf{N}}
\newcommand{\Sd}{\operatorname{Sd}}

\newcommand{\pol}{\operatorname{pol}}
\newcommand{\Id}{\operatorname{Id}}
\newcommand{\Ad}{\operatorname{Ad}}
\newcommand{\id}{\operatorname{id}}
\newcommand{\Tr}{\operatorname{Tr}}

\newcommand{\Aut}{\operatorname{Aut}}
\newcommand{\SL}{\operatorname{SL}}

\newcommand{\op}{\operatorname{op}}
\newcommand{\co}{\operatorname{co}}
\newcommand{\esssup}{\mathop{\mbox{ess sup}}}
\newcommand{\NC}{\operatorname{NC}}
\newcommand{\CC}{\operatorname{C}}
\newcommand{\Inn}{\operatorname{Inn}}
\newcommand{\Out}{\operatorname{Out}}
\newcommand{\fin}{\operatorname{fin}}
\newcommand{\cb}{\operatorname{cb}}
\newcommand{\fois}{\operatorname{times}}

\begin{document}

\title[Structural results for Free Araki-Woods Factors]{Structural results for free Araki-Woods factors and their continuous cores}

\begin{abstract}
We show that for any type ${\rm III_1}$ free Araki-Woods factor $\mathcal{M} = \Gamma(H_\R, U_t)''$ associated with an orthogonal representation $(U_t)$ of $\R$ on a separable real Hilbert space $H_\R$, the continuous core $M = \mathcal{M} \rtimes_\sigma \R$ is a semisolid ${\rm II_\infty}$ factor, i.e. for any non-zero finite projection $q \in M$, the ${\rm II_1}$ factor $qMq$ is semisolid. If the representation $(U_t)$ is moreover assumed to be mixing, then we prove that the core $M$ is solid. As an application, we construct an example of a non-amenable solid ${\rm II_1}$ factor $N$ with full fundamental group, i.e. $\mathcal{F}(N) = \R^*_+$, which is not isomorphic to any interpolated free group factor $L(\F_t)$, for $1 < t \leq +\infty$.
\end{abstract}

\author{Cyril Houdayer }

\address{CNRS-ENS Lyon \\
UMPA UMR 5669 \\
69364 Lyon cedex 7 \\
France}

\email{cyril.houdayer@umpa.ens-lyon.fr}

\subjclass[2000]{46L10; 46L54}

\keywords{Free Araki-Woods factors; (Semi)solid factors; Deformation/Rigidity techniques; Intertwining techniques; Spectral gap property}

\maketitle

\section{Introduction and statement of the main results}

The {\it free Araki-Woods factors} were introduced by Shlyakhtenko in \cite{shlya97}. In the context of {\it free probability} theory, these factors can be regarded as the analogs of the {\it hyperfinite} factors coming from the CAR\footnote{Canonical Anticommutation Relations} functor. To each real separable Hilbert space $H_\R$ together with an orthogonal representation $(U_t)$ of $\R$ on $H_\R$, one can associate a von Neumann algebra denoted by $\Gamma(H_\R, U_t)''$, called the {\it free Araki-Woods} von Neumann algebra. The von Neumann algebra $\Gamma(H_\R, U_t)''$ comes equipped with a unique {\it free quasi-free state} denoted by $\varphi_U$, which is always normal and faithful on $\Gamma(H_\R, U_t)''$ (see Section \ref{preliminaries} for a more detailed construction). If $\dim H_\R = 1$, then $\Gamma(\R, \Id)'' \cong L^\infty[0, 1]$. If $\dim H_\R \geq 2$, then $\mathcal{M} = \Gamma(H_\R, U_t)''$ is a full factor. In particular, $\mathcal{M}$ can never be of type ${\rm III_0}$. The type classification of these factors is the following:
\begin{enumerate}
\item $\mathcal{M}$ is a type ${\rm II_1}$ factor iff the representation $(U_t)$ is trivial: in that case the functor $\Gamma$ is Voiculescu's free Gaussian functor \cite{voiculescu92}. Then $\Gamma(H_\R, \Id)'' \cong L(\F_{\dim H_\R})$.
\item $\mathcal{M}$ is a type ${\rm III_\lambda}$ factor, for $0 < \lambda < 1$, iff the representation $(U_t)$ is $\frac{2\pi}{|\log \lambda|}$-periodic. 
\item $\mathcal{M}$ is a type ${\rm III_1}$ factor iff $(U_t)$ is non-periodic and non-trivial.
\end{enumerate}

Using free probability techniques, Shlyakhtenko obtained several remarkable classification results for $\Gamma(H_\R, U_t)''$. For instance, if the orthogonal representations $(U_t)$ are {\it almost periodic}, then the free Araki-Woods factors $\mathcal{M}  = \Gamma(H_\R, U_t)''$ are completely classified up to state-preserving $\ast$-isomorphism \cite{shlya97}: they only depend on Connes' invariant $\Sd(\mathcal{M})$ which is equal in that case to the (countable) subgroup $S_U \subset \R^*_+$ generated by the eigenvalues of $(U_t)$. Moreover, the {\it discrete} core $\mathcal{M} \rtimes_\sigma \widehat{S_U}$ (where $\widehat{S_U}$ is the compact group dual of $S_U$) is $\ast$-isomorphic to $L(\F_\infty) \bar{\otimes} \mathbf{B}(\ell^2)$. Shlyakhtenko showed in \cite{shlya98} that if $(U_t)$ is the left regular representation, then the {\it continuous} core $M = \mathcal{M} \rtimes_\sigma \R$ is isomorphic to $L(\F_\infty) \bar{\otimes} \mathbf{B}(\ell^2)$ and the dual {\textquotedblleft trace-scaling\textquotedblright} action $(\theta_s)$ is precisely the one constructed by R\u{a}dulescu \cite{radulescu1991}. For more on free Araki-Woods factors, we refer to \cite{{houdayer5}, {houdayer2}, {shlya2004b}, {shlya2004}, {shlya2003}, {shlya99}, {shlya98}, {shlya97}} and also to Vaes' Bourbaki seminar \cite{vaes2004}.

The free Araki-Woods factors as well as their continuous cores carry a {\it malleable deformation} in the sense of Popa. Then we will use the {\it deformation/rigidity} strategy together with the {\it intertwining techniques} in order to study the associated continuous cores. The high flexibility of this approach will allow us to work in a {\it semifinite} setting, so that we can obtain new structural/indecomposability results for the continuous cores of the free Araki-Woods factors. We first need to recall a few concepts. Following Ozawa \cite{{ozawa2003}, {ozawa2004}}, a finite von Neumann algebra $N$ is said to be:
\begin{itemize}
\item {\it solid} if for any diffuse von Neumann subalgebra $A \subset N$, the relative commutant $A' \cap N$ is amenable;
\item {\it semisolid} if for any type ${\rm II_1}$ von Neumann subalgebra $A \subset N$, the relative commutant $A' \cap N$ is amenable.
\end{itemize}
It is easy to check that solidity and semisolidity for ${\rm II_1}$ factors are stable under taking amplification by any $t > 0$. Moreover, if $N$ is a non-amenable ${\rm II_1}$ factor, then solid $\Longrightarrow$ semisolid $\Longrightarrow$ prime. Recall in this respect that $N$ is said to be {\it prime} if it cannot be written as the tensor product of two diffuse factors. 

Ozawa discovered a class $\mathcal{S}$ of countable groups for which whenever $\Gamma \in \mathcal{S}$, the group von Neumann algebra $L(\Gamma)$ is solid \cite{ozawa2003}. He showed that the following countable groups belong to the class $\mathcal{S}$: the word-hyperbolic groups \cite{ozawa2003}, the wreath products $\Lambda \wr \Gamma$ for $\Lambda$ amenable and $\Gamma \in \mathcal{S}$ \cite{ozawa2004}, and $\Z^2 \rtimes \SL(2, \Z)$ \cite{ozawa2008}. He moreover proved that if $\Gamma \in \mathcal{S}$, then for any free, ergodic, p.m.p. action $\Gamma \curvearrowright (X, \mu)$, the corresponding ${\rm II_1}$ factor $L^\infty(X, \mu) \rtimes \Gamma$ is semisolid \cite{ozawa2004}. Recall that a non-amenable solid ${\rm II_1}$ factor does not have property $\Gamma$ of Murray \& von Neumann \cite{ozawa2003}. 

\begin{df}
Let $M$ be a ${\rm II_\infty}$ factor and let $\Tr$ be a fixed faithful normal semifinite trace on $M$. We shall say that $M$ is {\it solid} (resp. {\it semisolid}) if for any non-zero projection $q \in M$ such that $\Tr(q) < \infty$, the ${\rm II_1}$ factor $qMq$ is solid (resp. semisolid).
\end{df}

Recall that an orthogonal/unitary representation $(U_t)$ acting on $H$ is said to be {\it mixing} if for any $\xi, \eta \in H$, $\langle U_t \xi, \eta\rangle \to 0$, as $|t| \to \infty$. The main result of this paper is the following:

\begin{theo}\label{mainresult}
Let $\mathcal{M} = \Gamma(H_\R, U_t)''$ be a type ${\rm III_1}$ free Araki-Woods factor. Then the continuous core $M = \mathcal{M} \rtimes_\sigma \R$ is a semisolid ${\rm II_\infty}$ factor. Since $M$ is non-amenable, $M$ is always a prime factor. If the representation $(U_t)$ is moreover assumed to be mixing, then $M$ is a solid ${\rm II_\infty}$ factor.
\end{theo}

The proof of Theorem \ref{mainresult} follows Popa's {\it
deformation/rigidity} strategy. This theory has been successfully used over the last eight years to give a plethora of new classification/rigidity results for crossed products/free products von Neumann algebras. We refer to \cite{{houdayer4}, {houdayer3}, {ipp}, {popasup}, {popasolid}, {popamal1}, {popa2001}, {popamsri}, {popavaes}, {vaesbern}} for some applications of the deformation/rigidity technique. We point out that in the present paper, the rigidity part does not rely on the notion of (relative) property (T) but rather on a certain {\it spectral gap} property discovered by Popa in \cite{{popasup}, {popasolid}}. Using this powerful technique, Popa was able to show for instance that the Bernoulli action of groups of the form $\Gamma_1 \times \Gamma_2$, with $\Gamma_1$ non-amenable and $\Gamma_2$ infinite is $\mathcal{U}_{\fin}$-cocycle superrigid \cite{popasup}. The spectral gap rigidity principle gave also a new approach to proving primeness and (semi)solidity for type ${\rm II_1}/{\rm III}$ factors \cite{{CI}, {houdayer4}, {popasup}, {popasolid}}. We briefly remind below the concepts that we will play against each other in order to prove Theorem \ref{mainresult}:
\begin{enumerate}
\item  The first ingredient we will use is the {\textquotedblleft malleable deformation\textquotedblright} by automorphisms $(\alpha_t, \beta)$ defined on the free Araki-Woods factor $\mathcal{M} \ast \mathcal{M} = \Gamma(H_\R \oplus H_\R, U_t \oplus U_t)''$. This deformation naturally arises as the {\textquotedblleft second quantization\textquotedblright} of the rotations/reflection defined on $H_\R \oplus H_\R$ that commute with $U_t \oplus U_t$. It was shown in \cite{popasup} that such a deformation automatically features a certain {\textquotedblleft transversality property\textquotedblright} (see Lemma $2.1$ in \cite{popasup}) which will be of essential use in our proof.

\item The second ingredient we will use is the {\em spectral gap rigidity} principle discovered by Popa in \cite{{popasup}, {popasolid}}. Let $B \subset M_i$ be an inclusion of finite von Neumann algebras, for $i = 1, 2$, with $B$ amenable. Write $M = M_1 \ast_B M_2$. Then for any von Neumann subalgebra $Q \subset M_1$ with no amenable direct summand, the action by conjugation $\Ad(\mathcal{U}(Q)) \curvearrowright M$ has {\textquotedblleft spectral gap\textquotedblright} relative to $M_1$: for any $\varepsilon > 0$, there exist $\delta > 0$ and a finite {\textquotedblleft critical\textquotedblright}  subset $F \subset \mathcal{U}(Q)$ such that for any $x \in (M)_1$ (the unit ball of $M$), if $\|u x u^* - x\|_2 \leq \delta$, $\forall u \in F$, then $\|x - E_{M_1}(x)\|_2 \leq \varepsilon$.

\item Let $M = \mathcal{M} \rtimes_\sigma \R$ be the continuous core the free Araki-Woods factor $\mathcal{M}$. Let $q \in M$ be a non-zero finite projection. A combination of $(1)$ and $(2)$ yields that for any $Q \subset qMq$ with no amenable direct summand, the malleable deformation $(\alpha_t)$ necessarily converges uniformly in $\left\| \cdot \right\|_2$ on $(Q' \cap qMq)_1$. Then, using Popa's intertwining techniques, one can locate the position of $Q' \cap qMq$ inside $qMq$.
\end{enumerate}

The second result of this paper provides a new example of a non-amenable solid ${\rm II_1}$ factor. We first need the following:

\begin{exam}\label{examplerep}
Using results of \cite{antoniou}, we construct an example of an orthogonal representation $(U_t)$ of $\R$ on a (separable) real Hilbert space $K_\R$ such that:
\begin{enumerate}
\item $(U_t)$ is mixing.
\item The spectral measure of $\bigoplus_{n \geq 1} U_t^{\otimes n}$ is singular w.r.t. the Lebesgue measure on $\R$.
\end{enumerate}
\end{exam}

Shlyakhtenko showed in \cite{shlya98} that if the spectral measure of the representation $\bigoplus_{n \geq 1} U_t^{\otimes n}$ is singular w.r.t. the Lebesgue measure, then the continuous core of the free Araki-Woods factor $\Gamma(H_\R, U_t)''$ cannot be isomorphic to any $L(\F_t) \bar{\otimes} \mathbf{B}(\ell^2)$, for $1 < t \leq \infty$, where $L(\F_t)$ denote the interpolated free group factors \cite{{dykema94}, {radulescu1994}}. Therefore, we obtain:

\begin{theo}\label{examplesolid}
Let $(U_t)$ be an orthogonal representation acting on $K_\R$ as in Example $\ref{examplerep}$. Denote by $\mathcal{M} = \Gamma(K_\R, U_t)''$ the corresponding free Araki-Woods factor and by $M = \mathcal{M}Ê\rtimes_\sigma \R$ its continuous core. Let $q \in L(\R)$ be a non-zero projection such that $\Tr(q) < \infty$. Then the non-amenable ${\rm II_1}$ factor $qMq$ is solid, has full fundamental group, i.e. $\mathcal{F}(qMq) = \R^*_+$, and is not isomorphic to any interpolated free group factor $L(\F_t)$, for $1 < t \leq \infty$.
\end{theo}

The paper is organized as follows. In Section \ref{preliminaries}, we recall the necessary background on free Araki-Woods factors as well as intertwining techniques for (semi)finite von Neumann algebras. Section \ref{results} is mainly devoted to the proof of Theorem \ref{mainresult}, following the deformation/spectral gap rigidity strategy presented above. In the last Section, we construct Example \ref{examplerep} and deduce Theorem \ref{examplesolid}.

{\bf Acknowledgement.} The author would like to thank Prof. D. Shlyakhtenko for suggesting him the idea of looking at indecomposability results for the continuous cores of the free Araki-Woods factors and also for the useful discussions.

\section{Preliminary background}\label{preliminaries}

\subsection{Shlyakhtenko's free Araki-Woods factors}

Let $H_{\R}$ be a real separable Hilbert space and let $(U_t)$ be an orthogonal representation of $\R$ on $H_{\R}$ such that the map $t \mapsto U_t$ is strongly continuous. Let $H_\C = H_{\R} \otimes_{\R} \C$ be the complexified Hilbert space. We shall still denote by $(U_t)$ the corresponding unitary representation of $\R$ on $H_\C$. Let $A$ be the infinitesimal generator of $(U_t)$ on $H_\C$ (Stone's theorem), so that $A$ is the positive, self-adjoint, (possibly) unbounded operator on $H_\C$ which satisfies $U_t = A^{it}$, for every $t \in \R$. Define another inner product on $H_\C$ by
\begin{equation*}
\langle \xi, \eta\rangle_{U} = \left\langle \frac{2}{1 + A^{-1}}\xi, \eta\right\rangle, \forall \xi, \eta \in H_\C.
\end{equation*}
Note that for any $\xi \in H_\R$, $\|\xi\|_U = \|\xi\|$; also, for any $\xi , \eta \in H_\R$, $\Re(\langle \xi, \eta\rangle_U) = \langle \xi, \eta\rangle$, where $\Re$ denotes the real part. Denote by $H$ the completion of $H_\C$ w.r.t. the new inner product $\langle \cdot, \cdot\rangle_U$, and note that $(U_t)$ is still a unitary representation on $H$. Introduce now the \emph{full Fock space} of $H$:
\begin{equation*}
\mathcal{F}(H) =\C\Omega \oplus \bigoplus_{n = 1}^{\infty} H^{\bar{\otimes} n}.
\end{equation*}
The unit vector $\Omega$ is called the \emph{vacuum vector}. For any $\xi \in H$, we have the \emph{left creation operator}
\begin{equation*}
\ell(\xi) : \mathcal{F}(H) \to \mathcal{F}(H) : \left\{ 
{\begin{array}{l} \ell(\xi)\Omega = \xi, \\ \ell(\xi)(\xi_1 \otimes \cdots \otimes \xi_n) = \xi \otimes \xi_1 \otimes \cdots \otimes \xi_n.
\end{array}} \right.
\end{equation*}
For any $\xi \in H$, we denote by $s(\xi)$ the real part of $\ell(\xi)$ given by
\begin{equation*}
s(\xi) = \frac{\ell(\xi) + \ell(\xi)^*}{2}.
\end{equation*}
The crucial result of Voiculescu \cite{voiculescu92} is that the distribution of the operator $s(\xi)$ w.r.t. the vacuum vector state $\varphi_U = \langle \cdot \Omega, \Omega\rangle_U$ is the semicircular law of Wigner supported on the interval $[-\|\xi\|, \|\xi\|]$.

\begin{df}[Shlyakhtenko, \cite{shlya97}]
Let $(U_t)$ be an orthogonal representation of $\R$ on the real Hilbert space $H_{\R}$. The \emph{free Araki-Woods} von Neumann algebra associated with $H_\R$ and $(U_t)$, denoted by $\Gamma(H_{\R}, U_t)''$, is defined by
\begin{equation*}
\Gamma(H_{\R}, U_t)'' := \{s(\xi) : \xi \in H_{\R}\}''.
\end{equation*}
The vector state $\varphi_{U} = \langle \cdot \Omega, \Omega\rangle_U$ is called the \emph{free quasi-free state}. It is normal and faithful on $\Gamma(H_\R, U_t)''$.
\end{df}

Recall that for any type ${\rm III_1}$ factor $\mathcal{M}$, Connes-Takesaki's continuous decomposition \cite{{connes73}, {takesaki73}} yields
\begin{equation*}
\mathcal{M} \bar{\otimes} \mathbf{B}(L^2(\R)) \cong (\mathcal{M} \rtimes_\sigma \R) \rtimes_\theta \R,
\end{equation*}
where the {\it continuous core} $\mathcal{M} \rtimes_\sigma \R$ is a ${\rm II_\infty}$ factor and $\theta$ is the {\it trace-scaling} action \cite{takesaki73}:
\begin{equation*}
\Tr(\theta_s(x)) = e^{-s} \Tr(x), \forall x \in (\mathcal{M} \rtimes_\sigma \R)_+, \forall s \in \R.
\end{equation*}
The fact that $\mathcal{M} \rtimes_\sigma \R$ does not depend on the choice of a f.n. state on $\mathcal{M}$ follows from Connes' Radon-Nikodym derivative theorem \cite{connes73}. Moreover, for any non-zero finite projection $q \in M = \mathcal{M} \rtimes_\sigma \R$, the ${\rm II_1}$ factor $qMq$ has full fundamental group.

Following \cite{connes74}, a factor $\mathcal{M}$ (with separable predual) is said to be {\it full} if the subgroup of inner automorphisms $\Inn(\mathcal{M}) \subset \Aut(\mathcal{M})$ is closed. Recall that $\Aut(\mathcal{M})$ is endowed with the $u$-topology: for any sequence $(\theta_n)$ in $\Aut(\mathcal{M})$, 
\begin{equation*}
\theta_n \to \Id, \mbox{ as } n \to \infty \Longleftrightarrow \left\| \varphi \circ \theta_n - \varphi \right\| \to 0, \mbox{ as } n \to \infty, \forall \varphi \in \mathcal{M}_*.
\end{equation*}
Since $\mathcal{M}$ has a separable predual, $\Aut(\mathcal{M})$ is a polish group. For any ${\rm II_1}$ factor $N$, $N$ is full iff $N$ does not have property $\Gamma$ of Murray \& von Neumann (see \cite{connes74}).

Denote by $\pi : \Aut(\mathcal{M}) \to \Out(\mathcal{M})$ the canonical projection. Assume $\mathcal{M}$ is a full factor so that $\Out(\mathcal{M})$ is a Hausdorff topological group. Fix a f.n. state $\varphi$ on $\mathcal{M}$. Connes' invariant $\tau(\mathcal{M})$ is defined as the {\it weakest topology} on $\R$ that makes the map
\begin{eqnarray*}
\R & \to & \Out(\mathcal{M}) \\
t & \mapsto & \pi\left( \sigma_t^\varphi \right)
\end{eqnarray*}
continuous. Note that this map does not depend on the choice of the f.n. state $\varphi$ on $\mathcal{M}$ \cite{connes73}.

Denote by $\mathcal F(U_t) = \bigoplus_{n \in \N} U_t^{\otimes n}$. The modular group $\sigma^{\varphi_U}$ of the free quasi-free state is given by: $\sigma_t^{\varphi_U} = \Ad(\mathcal F(U_{-t}))$, for any $t \in \R$. The free Araki-Woods factors provided many new examples of full factors of type {\rm III} \cite{{barnett95}, {connes73}, {shlya2004}}. We can summarize their general properties in the following theorem (see also Vaes' Bourbaki seminar \cite{vaes2004}):

\begin{theo}[Shlyakhtenko, \cite{{shlya2004}, {shlya99}, {shlya98}, {shlya97}}]\label{faw}
Let $(U_t)$ be an orthogonal representation of $\R$ on the real Hilbert space $H_{\R}$ with $\dim H_{\R} \geq 2$. Denote by $\mathcal{M} = \Gamma(H_{\R}, U_t)''$.
\begin{enumerate}
\item $\mathcal M$ is a full factor and Connes' invariant $\tau(\mathcal M)$ is the weakest topology on $\R$ that makes the map $t \mapsto U_t$ strongly continuous.
\item $\mathcal{M}$ is of type ${\rm II_1}$ iff $U_t = \id$ for every $t \in \R$. In this case, $\mathcal{M} \cong L(\F_{\dim(H_\R)})$.
\item $\mathcal{M}$ is of type ${\rm III_{\lambda}}$ $(0 < \lambda < 1)$ iff $(U_t)$ is periodic of period $\frac{2\pi}{|\log \lambda|}$.
\item $\mathcal{M}$ is of type ${\rm III_1}$ in the other cases.
\item $\mathcal{M}$ has almost periodic states iff $(U_t)$ is almost periodic.
\end{enumerate}
\end{theo}

Moreover, it follows from \cite{shlya2004b} that any free Araki-Woods factor $\mathcal{M}$ is {\it generalized solid} in the sense of \cite{VV}: for any diffuse von Neumann subalgebra $A \subset \mathcal M$ such that there exists a faithful normal conditional expectation $E : \mathcal M \to A$, the relative commutant $A' \cap \mathcal M$ is amenable.

Notice that the centralizer of the free quasi-free state $\mathcal{M}^{\varphi_U}$ may be trivial. This is the case for instance when the representation $(U_t)$ has no eigenvectors. Nevertheless, the author recently proved in \cite{houdayer5} that for any type ${\rm III_1}$ free Araki-Woods factor $\mathcal{M}$, the {\it bicentralizer} is trivial, i.e. there always exists a faithful normal state $\psi$ on $\mathcal{M}$ such that $(\mathcal{M}^\psi)' \cap \mathcal{M} = \C$. We refer to \cite{haagerup84} for more on Connes' bicentralizer problem.

\begin{rem}[\cite{shlya97}]
Explicitly the value of $\varphi_U$ on a word in $s(\xi_\iota)$ is given by
\begin{equation}\label{formula}
\varphi_U(s(\xi_1) \cdots s(\xi_n)) = 2^{-n}\sum_
{(\{\beta_i, \gamma_i\}) 
 \in \NC(n), 
\beta_i < \gamma_i}
\prod_{k = 1}^{n/2}\langle \xi_{\beta_k}, \xi_{\gamma_k}\rangle_U.
\end{equation}
for $n$ even and is zero otherwise. Here $\NC(2p)$ stands for all the non-crossing pairings of the set $\{1, \dots, 2p\}$, i.e. pairings for which whenever $a < b < c < d$, and $a, c$ are in the same class, then $b, d$ are not in the same class. The total number of such pairings is given by the $p$-th Catalan number
\begin{equation*}
C_p = \frac{1}{p + 1}\begin{pmatrix}
2p \\
p
\end{pmatrix}.
\end{equation*}
\end{rem}

Recall that a continuous $\varphi$-preserving action $(\sigma_t)$ of $\R$ on a von Neumann algebra $\mathcal{M}$ endowed with a f.n. state $\varphi$ is said to be $\varphi$-{\it mixing} if for any $x, y \in \mathcal{M}$ with $\varphi(x) = \varphi(y) = 0$,
\begin{equation}\label{mixing}
\varphi(\sigma_t(x) y) \to 0, \mbox{ as } |t| \to \infty.
\end{equation}

\begin{prop}
Let $\mathcal{M} = \Gamma(H_\R, U_t)''$ be any free Araki-Woods factor and let $\varphi_U$ be the free quasi-free state. Then
\begin{equation*}
(U_t) \mbox{ is mixing } \Longleftrightarrow (\sigma_t^{\varphi_U}) \mbox{ is } \varphi_U\mbox{-mixing}.
\end{equation*}
\end{prop}

\begin{proof} We prove both directions.

$\Longleftarrow$ For any $\xi, \eta \in H_\R$, $\varphi_U(s(\xi)) = \varphi_U(s(\eta)) = 0$. Moreover,
\begin{eqnarray*}
\langle U_t \xi, \eta \rangle_U & = & 4 \varphi_U(s(U_t \xi) s(\eta)) \\
& = & 4 \varphi_U( \sigma_{-t}^{\varphi_U}(s(\xi)) s(\eta)) \to 0, \mbox{ as } |t| \to \infty. 
\end{eqnarray*}
It follows that $(U_t)$ is mixing.

$\Longrightarrow$ One needs to show that for any $x, y \in \mathcal{M}$,
\begin{equation*}
\lim_{|t| \to \infty} \varphi_U(\sigma_t^{\varphi_U}(x)y) = \varphi_U(x) \varphi_U(y).
\end{equation*}
Note that 
\begin{equation*}
\mbox{span} \left\{ 1, s(\xi_1) \cdots s(\xi_n) : n \geq 1, \xi_1, \dots, \xi_n \in H_\R \right\}
\end{equation*}
is a unital $\ast$-strongly dense $\ast$-subalgebra of $\mathcal{M}$. Using Kaplansky density theorem, it suffices to check Equation $(\ref{mixing})$ for $x, y \in \mathcal{M}$ of the following form:
\begin{eqnarray*}
x & = & s(\xi_1) \cdots s(\xi_{p}) \\
y & = & s(\eta_1) \cdots s(\eta_{q}).
\end{eqnarray*}

Assume that $p + q$ is odd. Then $p$ or $q$ is odd and we have $\varphi_U(\sigma_t^{\varphi_U}(x)y) = 0 = \varphi_U(x) \varphi_U(y)$, for any $t \in \R$.

Assume now that $p + q$ is even.
\begin{enumerate}
\item Suppose that $p, q$ are odd and write $p = 2k + 1$, $q = 2l + 1$. Then
\begin{eqnarray*}
\varphi_U(\sigma_t^{\varphi_U}(x) y) & = & \varphi_U(s(U_{-t} \xi_1) \cdots s(U_{-t} \xi_{2k + 1}) s(\eta_1) \cdots s(\eta_{2l + 1})) \\
& = & 2^{-2(k + l + 1)}\sum_
{(\{\beta_i, \gamma_i\}) 
 \in \NC(2(k+ l + 1)), 
\beta_i < \gamma_i}
\prod_{j = 1}^{k+ l+ 1}\langle h_{\beta_j}, h_{\gamma_j}\rangle_U,
\end{eqnarray*}
where the letter $h$ stands for $U_{-t} \xi$ or $\eta$. Notice that since $2k + 1$ and $2l + 1$ are odd, for any non-crossing pairing $(\{\beta_i, \gamma_i\}) 
 \in \NC(2(k+ l + 1))$, there must exist some $j \in \{1, \dots, k + l + 1\}$ such that $\langle h_{\beta_j}, h_{\gamma_j}\rangle = \langle U_{-t} \xi_{\beta_j}, \eta_{\gamma_j} \rangle$. Since we assumed that $(U_t)$ is mixing, it follows that $\varphi_U(\sigma_t^{\varphi_U}(x) y) \to 0 = \varphi_U(x) \varphi_U(y)$, as $|t| \to \infty$.
 
\item Suppose that $p, q$ are even and write $p = 2k$, $q = 2l$. Then
\begin{eqnarray*}
\varphi_U(\sigma_t^{\varphi_U}(x) y) & = & \varphi_U(s(U_{-t} \xi_1) \cdots s(U_{-t} \xi_{2k}) s(\eta_1) \cdots s(\eta_{2l})) \\
& = & 2^{-2(k + l)}\sum_
{(\{\beta_i, \gamma_i\}) 
 \in \NC(2(k+ l)), 
\beta_i < \gamma_i}
\prod_{j = 1}^{k+ l}\langle h_{\beta_j}, h_{\gamma_j}\rangle_U,
\end{eqnarray*}
where the letter $h$ stands for $U_{-t} \xi$ or $\eta$. Note that for a non-crossing pairing $\nu = (\{\beta_i, \gamma_i\}) \in \NC(2(k+ l))$ such that an element of $\{1, \dots, 2k\}$ and an element of $\{1, \dots, 2l\}$ are in the same class, the proof of $(1)$ yields that the corresponding product $\prod_{j = 1}^{k+ l}\langle h_{\beta_j}, h_{\gamma_j}\rangle_U$ goes to $0$, as $|t| \to \infty$. Thus, we just need to sum up over the non-crossing pairings $\nu$ of the form $\nu_1 \times \nu_2$, where $\nu_1$ is a non-crossing pairing on the set $\{1, \dots, 2k\}$ and $\nu_2$ is a non-crossing pairing on the set $\{1, \dots, 2l\}$. Consequently, we get $\varphi_U(\sigma_t^{\varphi_U}(x)y) \to \varphi_U(x) \varphi_U(y)$, as $|t| \to \infty$. 
\end{enumerate}
Therefore, $(\sigma_t^{\varphi_U})$ is mixing.
\end{proof}

\begin{prop}
Let $\mathcal{M} = \Gamma(H_\R, U_t)''$. If $(U_t)$ is mixing, then Connes' invariant $\tau(\mathcal{M})$ is the usual topology on $\R$.
\end{prop}

\begin{proof}
Let $\mathcal{M} = \Gamma(H_\R, U_t)''$. Recall from Theorem \ref{faw} that $\tau(\mathcal{M})$ is the weakest topology on $\R$ that makes the map $t \mapsto U_t$ strongly continuous. Let $(t_k)$ be a sequence in $\R$ such that $t_k \to 0$ w.r.t. the topology $\tau(\mathcal{M})$, as $k \to \infty$, i.e. $U_{t_k} \to \Id$ strongly, as $k \to \infty$. Fix $\xi \in H_\R$, $\|\xi\| = 1$. Since 
\begin{equation*}
\lim_{k \to \infty} \langle U_{t_k} \xi, \xi\rangle = 1
\end{equation*}
and $(U_t)$ is assumed to be mixing, it follows that $(t_k)$ is necessarily bounded. Let $t \in \R$ be any cluster point for the sequence $(t_k)$. Then $U_t = \Id$. Since $(U_t)$ is mixing, it follows that $t = 0$. Therefore $(t_k)$ converges to $0$ w.r.t. the usual topology on $\R$.
\end{proof}

\subsection{Intertwining techniques for (semi)finite von Neumann algebras}

Let $(B, \tau)$ be a finite von Neumann algebra with a distinguished f.n. trace. Since $\tau$ is fixed, we simply denote $L^2(B, \tau)$ by $L^2(B)$. Let $H$ be a right Hilbert $B$-module, i.e. $H$ is a complex (separable) Hilbert space together with a normal $\ast$-representation $\pi : B^{\op} \to \mathbf{B}(H)$. For any $b \in B$, and $\xi \in H$, we shall simply write $\pi(b^{\op}) \xi = \xi b$. By the general theory, we know that there exists an isometry $v : H \to \ell^2 \bar{\otimes} L^2(B)$ such that $v(\xi b) = v(\xi) b$, for any $\xi \in H$, $b \in B$. Since $p = vv^*$ commutes with the right $B$-action on $\ell^2 \bar{\otimes} L^2(B)$, it follows that $p \in \mathbf{B}(\ell^2) \bar{\otimes} B$. Thus, as right $B$-modules, we have $H_B \simeq p(\ell^2 \bar{\otimes} L^2(B))_B$.

On $\mathbf{B}(\ell^2) \bar{\otimes} B$, we define the following f.n. semifinite trace $\Tr$ (which depends on $\tau$): for any $x = [x_{ij}]_{i, j} \in (\mathbf{B}(\ell^2) \bar{\otimes} B)_+$,
\begin{equation*}
\Tr \left( [x_{ij}]_{i, j} \right) = \sum_i \tau(x_{ii}).
\end{equation*}

We set $\dim(H_B) = \Tr(vv^*)$. Note that the dimension of $H$ depends on $\tau$ but does not depend on the isometry $v$. Indeed take another isometry $w : H \to \ell^2 \bar{\otimes} L^2(B)$, satisfying $w(\xi b) = w(\xi) b$, for any $\xi \in H$, $b \in B$. Note that $vw^* \in \mathbf{B}(\ell^2) \bar{\otimes} B$ and $w^*w = v^*v = 1$. Thus, we have
\begin{equation*}
\Tr(vv^*) = \Tr(v w^*w v^*) = \Tr(w v^* v w^*) = \Tr(ww^*).
\end{equation*}

Assume that $\dim(H_B) < \infty$. Then for any $\varepsilon > 0$, there exists a central projection $z \in \mathcal{Z}(B)$, with $\tau(z) \geq 1 - \varepsilon$, such that the right $B$-module $Hz$ is finitely generated, i.e. of the form $p L^2(B)^{\oplus n}$ for some projection $p \in \mathbf{M}_n(\C) \otimes B$. The non-normalized trace on $\mathbf{M}_n(\C)$ will be denoted by $\Tr_n$. For simplicity, we shall denote $B^n := \mathbf{M}_n(\C) \otimes B$.

In \cite{{popamal1}, {popa2001}}, Popa introduced a powerful tool to prove the unitary conjugacy of two von Neumann subalgebras of a tracial von Neumann algebra $(M, \tau)$. If $A, B \subset (M, \tau)$ are two (possibly non-unital) von Neumann subalgebras, denote by $1_A, 1_B$ the units of $A$ and $B$. Note that we endow the finite von Neumann algebra $B$ with the trace $\tau(1_B \cdot 1_B) / \tau(1_B)$.

\begin{theo}[Popa, \cite{{popamal1}, {popa2001}}]\label{intertwining1}
Let $A, B \subset (M, \tau)$ be two (possibly non-unital) embeddings. The following are equivalent:
\begin{enumerate}
\item There exist $n \geq 1$, a (possibly non-unital) $\ast$-homomorphism $\psi : A \to B^n$ and a non-zero partial isometry $v \in \mathbf{M}_{1, n}(\C) \otimes 1_AM1_B$ such that $x v = v \psi(x)$, for any $x \in A$.

\item The bimodule $\vphantom{}_AL^2(1_AM1_B)_B$ contains a non-zero sub-bimodule $\vphantom{}_AH_B$ which satisfies $\dim(H_B) < \infty$. 

\item There is no sequence of unitaries $(u_k)$ in $A$ such that $\left\|E_B(a^* u_k b)\right\|_2 \to 0$, as $k \to \infty$, for any $a, b \in 1_A M 1_B$.
\end{enumerate}
\end{theo}
If one of the previous equivalent conditions is satisfied, we shall say that $A$ {\it embeds into} $B$ {\it inside} $M$ and denote $A \preceq_M B$.

\begin{df}[Popa \& Vaes, \cite{popavaes}]\label{weakmixing}
Let $A \subset B \subset (N, \tau)$ be an inclusion of finite von Neumann algebras. We say that $B \subset N$ is {\it weakly mixing through} $A$ if there exists a sequence of unitaries $(u_k)$ in $A$ such that
\begin{equation*}
\left\| E_B(a^* u_k b) \right\|_2 \to 0, \mbox { as } k \to \infty, \forall a,b \in N \ominus B.
\end{equation*}
\end{df}

The following result will be a crucial tool in Section \ref{results}: it will allow us to control the relative commutant $A' \cap N$ of certain subalgebras $A$ of a given von Neumann algebra $N$.

\begin{theo}[Popa, \cite{popamal1}]\label{principle}
Let $(N, \tau)$ be a finite von Neumann algebra and $A \subset B \subset N$ be von Neumann subalgebras. Assume that $B \subset N$ is weakly mixing through $A$. Then for any sub-bimodule $\vphantom{}_AH_B$ of $\vphantom{}_AL^2(N)_B$ such that $\dim(H_B) < \infty$, one has $H \subset L^2(B)$. In particular, $A' \cap N \subset B$.
\end{theo}

For our purpose, we will need to use Popa's intertwining techniques for {\it semifinite} von Neumann algebras. We refer to Section $2$ of \cite{houdayer4} where such techniques were developed. Namely, let $(M, \Tr)$ be a von Neumann algebra endowed with a faithful normal semifinite trace $\Tr$. We shall simply denote by $L^2(M)$ the $M$-$M$ bimodule $L^2(M, \Tr)$, and by $\left\| \cdot \right\|_{2, \Tr}$ the $L^2$-norm associated with the trace $\Tr$. We will use quite often the following inequality:
\begin{equation*}
\left\| x \eta y \right\|_{2, \Tr} \leq \left\| x \right\|_\infty \left\| y \right\|_\infty \left\| \eta \right\|_{2, \Tr}, \forall \eta \in L^2(M), \forall x, y \in M,
\end{equation*}
where $\left\| \cdot \right\|_\infty$ denotes the operator norm. We shall say that a projection $p \in M$ is $\Tr$-{\it finite} if $\Tr(p) < \infty$. Note that a non-zero $\Tr$-finite projection $p$ is necessarily finite and $\Tr(p \cdot p)/\Tr(p)$ is a f.n. (finite) trace on $pMp$. Remind that for any projections $p, q \in M$, we have $p \vee q - p \sim q - p \wedge q$. Then it follows that for any $\Tr$-finite projections $p, q \in M$, $p \vee q$ is still $\Tr$-finite and $\Tr(p \vee q) = \Tr(p) + \Tr(q) - \Tr(p \wedge q)$.

Note that if a sequence $(x_k)$ in $M$ converges to $0$ strongly, as $k \to \infty$, then for any non-zero $\Tr$-finite projection $q \in M$, $\left\| x_k q \right\|_{2, \Tr} \to 0$, as $k \to \infty$. Indeed,
\begin{eqnarray*}
x_k \to 0 \; \mbox{strongly in } M & \Longleftrightarrow & x^*_kx_k \to 0 \mbox{ weakly in } M \\
& \Longrightarrow & qx^*_kx_kq \to 0 \mbox{ weakly in } qMq \\
& \Longrightarrow & \Tr(qx^*_kx_kq) \to 0 \\
& \Longrightarrow & \left\| x_k q \right\|_{2, \Tr} \to 0.
\end{eqnarray*}
Moreover, there always exists an increasing sequence of $\Tr$-finite projections $(p_k)$ in $M$ such that $p_k \to 1$ strongly, as $k \to \infty$.

\begin{theo}[\cite{houdayer4}]\label{intertwining}
Let $(M, \Tr)$ be a semifinite von Neumann algebra. Let $B \subset M$ be a von Neumann subalgebra such that $\Tr_{|B}$ is still semifinite. Denote by $E_B : M \to B$ the unique $\Tr$-preserving faithful normal conditional expectation. Let $q \in M$ be a non-zero $\Tr$-finite projection. Let $A \subset qMq$ be a von Neumann subalgebra. The following conditions are equivalent:
\begin{enumerate}
\item There exists a $\Tr$-finite projection $p \in B$, $p \neq 0$, such that the bimodule $\vphantom{}_AL^2(qMp)_{pBp}$ contains a non-zero sub-bimodule $\vphantom{}_AH_{pBp}$ which satisfies $\dim(H_{pBp}) < \infty$, where $pBp$ is endowed with the finite trace $\Tr(p \cdot p) / \Tr(p)$.

\item There is no sequence of unitaries $(u_k)$ in $A$ such that $E_B(x^* u_k y) \to 0$ strongly, as $k \to \infty$, for any $x, y \in qM$.
\end{enumerate}
\end{theo}

\begin{df}
Under the assumptions of Theorem $\ref{intertwining}$, if one of the equivalent conditions is satisfied, we shall still say that $A$ {\it embeds into} $B$ {\it inside} $M$ and still denote $A \preceq_M B$.
\end{df}

\section{Structural results for the continuous cores of $\Gamma(H_\R, U_t)''$}\label{results}

\subsection{Deformation/spectral gap rigidity strategy}

We first introduce some notation we will be using throughout this section. Let $H_\R$ be a separable real Hilbert space ($\dim(H_\R) \geq 2$) and let $(U_t)$ be an orthogonal representation of $\R$ on $H_\R$ that we assume to be neither trivial nor periodic. We set: 
\begin{itemize}
\item $\mathcal{M} = \Gamma(H_\R, U_t)''$ is the free Araki-Woods factor associated with $(H_\R, U_t)$,  $\varphi$ is the free quasi-free state and $\sigma$ is the modular group of the state $\varphi$.  $\mathcal{M}$ is necessarily a type ${\rm III_1}$ factor since $(U_t)$ is neither periodic nor trivial. 
\item $M = \mathcal{M} \rtimes_\sigma \R$ is the continuous core of $\mathcal{M}$ and $\Tr$ is the semifinite trace associated with the state $\varphi$. $M$ is a ${\rm II_\infty}$ factor since $\mathcal{M}$ is a type ${\rm III_1}$ factor.
\item Likewise $\widetilde{\mathcal{M}} = \Gamma(H_\R \oplus H_\R, U_t \oplus U_t)''$, $\widetilde{\varphi}$ is the corresponding free quasi-free state and $\widetilde{\sigma}$ is the modular group of $\widetilde{\varphi}$.
\item $\widetilde{M} = \widetilde{\mathcal{M}} \rtimes_{\widetilde{\sigma}} \R$ is the continuous core of $\widetilde{\mathcal{M}}$ and $\widetilde{\Tr}$ is the f.n. semifinite trace associated with $\widetilde{\varphi}$.
\end{itemize}
It follows from \cite{shlya97} that
\begin{equation*}
\left( \widetilde{\mathcal{M}}, \widetilde{\varphi} \right) \cong (\mathcal{M}, \varphi) \ast (\mathcal{M}, \varphi).
\end{equation*}
In the latter free product, we shall write $\mathcal{M}_1$ for the first copy of $\mathcal{M}$ and $\mathcal{M}_2$ for the second copy of $\mathcal{M}$. We regard $\mathcal{M} \subset \widetilde{\mathcal{M}}$ via the identification of $\mathcal{M}$ with $\mathcal{M}_1$.

Denote by $(\lambda_t)$ the unitaries in $L(\R)$ that implement the modular action $\sigma$ on $\mathcal{M}$ (resp. $\widetilde{\sigma}$ on $\widetilde{\mathcal{M}}$). Define the following faithful normal conditional expectations:
\begin{itemize}
\item $E : M \to L(\R)$ such that $E(x \lambda_t) = \varphi(x) \lambda_t$, for every $x \in \mathcal{M}$ and $t \in \R$;
\item $\widetilde{E} : \widetilde{M} \to L(\R)$ such that $\widetilde{E}(x \lambda_t) = \widetilde{\varphi}(x) \lambda_t$, for every $x \in \widetilde{\mathcal{M}}$ and $t \in \R$.
\end{itemize}
Then
\begin{equation*}
\left( \widetilde{M}, \widetilde{E} \right) \cong (M, E) \ast_{L(\R)} (M, E).
\end{equation*}
Likewise, in the latter amalgamated free product, we shall write $M_1$ for the first copy of $M$ and $M_2$ for the second copy of $M$. We regard $M \subset \widetilde{M}$ via the identification of $M$ with $M_1$. Notice that the conditional expectation $E$ (resp. $\widetilde{E}$) preserves the canonical semifinite trace $\Tr$ (resp. $\widetilde{\Tr}$) associated with the state $\varphi$ (resp. $\widetilde{\varphi}$) (see \cite{ueda}).

Consider the following orthogonal representation of $\R$ on $H_\R \oplus H_\R$:
\begin{equation*}
V_s = \begin{pmatrix}
\cos(\frac{\pi}{2}s) & -\sin(\frac{\pi}{2}s) \\
\sin(\frac{\pi}{2}s) & \cos(\frac{\pi}{2}s)
\end{pmatrix}, \forall s \in \R.
\end{equation*}
Let $(\alpha_s)$ be the natural action on $\left( \widetilde{\mathcal{M}}, \widetilde{\varphi} \right)$ associated with $(V_s)$: $\alpha_s = \Ad(\mathcal{F}(V_s))$, for every $s \in \R$. In particular, we have
\begin{equation*}
\alpha_s(s\begin{pmatrix}
\xi \\ \eta \end{pmatrix}) =  s(V_s\begin{pmatrix}
\xi \\ \eta \end{pmatrix}), \forall s \in \R, \forall \xi, \eta \in H_\R,
\end{equation*}
and the action $(\alpha_s)$ is $\widetilde{\varphi}$-preserving. We can easily see that the representation $(V_s)$ commutes with the representation $(U_t \oplus U_t)$. Consequently,  $(\alpha_s)$ commutes with modular action $\widetilde{\sigma}$. Moreover, $\alpha_1(x \ast 1) = 1 \ast x$, for every $a \in \mathcal{M}$. At last,  consider the automorphism $\beta$ defined on $\left( \widetilde{\mathcal{M}}, \widetilde{\varphi} \right)$ by:
\begin{equation*}
\beta(s\begin{pmatrix}
\xi \\ \eta \end{pmatrix}) = s\begin{pmatrix}
\xi \\ -\eta \end{pmatrix}, \forall \xi, \eta \in H_\R. 
\end{equation*}
It is straightforward to check that $\beta$ commutes with the modular action $\widetilde{\sigma}$, $\beta^2 = \Id$, $\beta_{|\mathcal{M}} = \Id_{|\mathcal{M}}$ and $\beta\alpha_{s} = \alpha_{-s}\beta$, $\forall s \in \R$. Since $(\alpha_s)$ and $\beta$ commute with the modular action $\widetilde{\sigma}$, one may extend $(\alpha_s)$ and $\beta$ to $\widetilde{M}$ by ${\alpha_s}_{|L(\R)} = \Id_{L(\R)}$, for every $s \in \R$ and $\beta_{|L(\R)} = \Id_{L(\R)}$. Moreover $(\alpha_s, \beta)$ preserves the semifinite trace $\widetilde{\Tr}$. Let's summarize what we have done so far:

\begin{prop}
The $\widetilde{\Tr}$-preserving deformation $(\alpha_s, \beta)$ defined on $\widetilde{M}$ is \emph{s-malleable}:
\begin{enumerate}
\item ${\alpha_s}_{|L(\R)} = \Id_{L(\R)}$, for every $s \in \R$ and $\alpha_1(x \ast_{L(\R)} 1) = 1 \ast_{L(\R)} x$, for every $x \in M$. 
\item $\beta^2 = \Id$ and $\beta_{|M} = \Id_{|M}$. 
\item $\beta \alpha_s = \alpha_{-s} \beta$, for every $s \in \R$.
\end{enumerate}
\end{prop}

Denote by $E_M : \widetilde{M} \to M$ the canonical trace-preserving conditional expectation. Since $\widetilde{\Tr}_{|M} = \Tr$, we will simply denote by $\Tr$ the semifinite trace on $\widetilde{M}$. 
Remind that the s-malleable deformation $(\alpha_s, \beta)$ automatically features a certain {\it transversality property}. 

\begin{prop}[Popa, \cite{popasup}]\label{transversality}
We have the following:
\begin{equation}\label{trans}
\left\| x - \alpha_{2s}(x) \right\|_{2, \Tr} \leq 2 \left\| \alpha_s(x) - E_{M}(\alpha_s(x)) \right\|_{2, \Tr}, \; \forall x \in L^2(M, \Tr), \forall s > 0.
\end{equation}
\end{prop}

The next proposition refered in the Introduction as the {\it spectral gap} property was first proved by Popa in \cite{popasolid} for free products of finite von Neumann algebras. We will need the following straightforward generalization: 

\begin{prop}[\cite{houdayer4}]\label{spectralgap}
We keep the same notation as before. Let $q \in M$ be a non-zero projection such that $\Tr(q) < \infty$. Let $Q \subset qMq$ be a von Neumann subalgebra with no amenable direct summand. Then for any free ultrafilter $\omega$ on $\N$, we have $Q' \cap (q\widetilde{M}q)^\omega \subset (qMq)^\omega$.
\end{prop}

Let $q \in M$ be a non-zero projection such that $\Tr(q) < \infty$. Note that $\Tr(q \cdot q)/\Tr(q)$ is a finite trace on $q \widetilde{M} q$. If $Q \subset qMq$ has no amenable direct summand, then for any $\varepsilon > 0$, there exist $\delta > 0$ and a finite subset $F \subset \mathcal{U}(Q)$ such that for any $x \in (q\widetilde{M}q)_1$ (the unit ball w.r.t. the operator norm),
\begin{equation}\label{spectral}
\left\| ux - xu \right\|_{2, \Tr} < \delta, \forall u \in F \Longrightarrow \left\| x - E_{qMq}(x) \right\|_{2, \Tr} < \varepsilon.
\end{equation}
We will simply denote $ux - xu$ by $[u, x]$.

\subsection{Semisolidity of the continuous core}

The following theorem is in some ways a reminiscence of a result of Ioana, Peterson \& Popa, namely Theorem $4.3$ of \cite{ipp} and also Theorem $4.2$ of \cite{houdayer4}. The deformation/spectral gap rigidity strategy enables us to locate inside the core $M$ of a free Araki-Woods factor the position of subalgebras $A \subset M$ with a {\it large} relative commutant $A' \cap M$.

\begin{theo}\label{semisolidity}
Let $\mathcal{M} = \Gamma(H_\R, U_t)''$ be a free Araki-Woods factor and $M = \mathcal{M} \rtimes_\sigma \R$ be its continuous core. Let $q \in L(\R) \subset M$ be a non-zero projection such that $\Tr(q) < \infty$. Let $Q \subset qMq$ be a von Neumann subalgebra with no amenable direct summand. Then $Q' \cap qMq \preceq_M L(\R)$.
\end{theo}

\begin{cor}\label{semisoliditybis}
Let $\mathcal{M} = \Gamma(H_\R, U_t)''$ be a free Araki-Woods factor of type ${\rm III_1}$. Then the continuous core $M = \mathcal{M} \rtimes_\sigma \R$ is a semisolid ${\rm II_\infty}$ factor. Since $M$ is non-amenable, $M$ is always a prime factor.
\end{cor}

\begin{proof}[Proof of Theorem $\ref{semisolidity}$]
Let $q \in L(\R)$ be a non-zero projection such that $\Tr(q) < \infty$. Let $Q \subset qMq$ be a von Neumann subalgebra with no amenable direct summand. Denote by $Q_0 = Q' \cap qMq$. We keep the notation introduced previously and regard $M \subset \widetilde{M} = M_1 \ast_{L(\R)} M_2$ via the identification of $M$ with $M_1$. Remind that ${\alpha_s}_{|L(\R)} = \Id_{L(\R)}$, for every $s \in \R$. In particular $\alpha_s(q) = q$, for every $s \in \R$.

{\bf Step (1) : Using the spectral gap condition and the transversality property of $(\alpha_t, \beta)$ to find $t > 0$ and a non-zero intertwiner $v$ between $\Id$ and $\alpha_t$.}

Let $\varepsilon = \frac{1}{4}\left\| q \right\|_{2, \Tr}$. We know that there exist $\delta > 0$ and a finite subset $F \subset \mathcal{U}(Q)$, such that for every $x \in (q \widetilde{M} q)_1$,
\begin{equation*}
\left\| [x, u] \right\|_{2, \Tr} \leq \delta, \forall u \in F \Longrightarrow \left\| x - E_{qMq}(x) \right\|_{2, \Tr} \leq \varepsilon.
\end{equation*}
Since $\alpha_t \to \Id$ pointwise $\ast$-strongly, as $t \to 0$, and since $F$ is a finite subset of $Q \subset qMq$, we may choose $t = 1/2^k$ small enough ($k \geq 1$) such that 
\begin{equation*}
\max\{\left\| u - \alpha_t(u) \right\|_{2, \Tr} : u \in F\} \leq \frac{\delta}{2}.
\end{equation*}
For every $x \in (Q_0)_1$ and every $u \in F \subset Q$, since $[u, x] = 0$, we have
\begin{eqnarray*}
\left\| [\alpha_t(x), u] \right\|_{2, \Tr} & = &  \left\| [\alpha_t(x), u - \alpha_t(u)] \right\|_{2, \Tr} \\
& \leq &  2\|u - \alpha_t(u)\|_{2, \Tr} \\
& \leq & \delta.
\end{eqnarray*}
Consequently, we get for every $x \in (Q_0)_1$, $\left\| \alpha_t(x) - E_{qMq}(\alpha_t(x)) \right\|_{2, \Tr} \leq \varepsilon$. Using Proposition $\ref{transversality}$, we obtain for every $x \in (Q_0)_1$
\begin{equation*}
\left\| x - \alpha_s(x) \right\|_{2, \Tr} \leq \frac{1}{2}\left\| q \right\|_{2, \Tr},
\end{equation*}
where $s = 2t$. Thus, for every $u \in \mathcal{U}(Q_0)$, we have
\begin{eqnarray*}
\left\| u^*\alpha_s(u) - q \right\|_{2, \Tr} & = & \left\| u^*(\alpha_s(u) - u) \right\|_{2, \Tr} \\
& \leq & \left\| u - \alpha_s(u) \right\|_{2, \Tr} \\
& \leq & \frac{1}{2} \left\| q \right\|_{2, \Tr}.
\end{eqnarray*}
Denote by $\mathcal{C} = \overline{\co}^w \{u^*\alpha_s(u) : u \in \mathcal{U}(Q_0)\} \subset q L^2(\widetilde{M})q$ the ultraweak closure of the convex hull of all $u^*\alpha_s(u)$, where $u \in \mathcal{U}(Q_0)$. Denote by $a$ the unique element in $\mathcal{C}$ of minimal $\left\| \cdot \right\|_{2, \Tr}$-norm. Since $\left\| a - q \right\|_{2, \Tr} \leq 1/2 \left\| q \right\|_{2, \Tr}$, necessarily $a \neq 0$. Fix $u \in \mathcal{U}(Q_0)$. Since $u^* a \alpha_s(u) \in \mathcal{C}$ and $\left\| u^* a \alpha_s(u) \right\|_{2, \Tr} = \left\| a \right\|_{2, \Tr}$, necessarily $u^* a \alpha_s(u) = a$. Taking $v = \pol(a)$ the polar part of $a$, we have found a non-zero partial isometry $v \in q\widetilde{M}q$ such that
\begin{equation}\label{specgap}
x v = v \alpha_s (x), \forall x \in Q_0.
\end{equation}

{\bf Step (2) : Proving $Q_0 \preceq_M L(\R)$ using the malleability of $(\alpha_t, \beta)$.} By contradiction, assume $Q_0 \npreceq_M L(\R)$. The first task is to lift Equation $(\ref{specgap})$ to $s = 1$. Note that it is enough to find a non-zero partial isometry $w \in q\widetilde{M}q$ such that
\begin{equation*}
x w = w \alpha_{2s} (x), \forall x \in Q_0.
\end{equation*}
Indeed, by induction we can go till $s = 1$ (because $s = 1/2^{k - 1}$). Remind that $\beta(z) = z$, for every $z \in M$. Note that $vv^* \in Q_0' \cap q\widetilde{M}q$. Since $Q_0 \npreceq_{M} L(\R)$, we know from Theorem $2.4$ in \cite{houdayer4} that $Q_0' \cap q\widetilde{M}q \subset qMq$. In particular, $vv^* \in qMq$. Set $w = \alpha_s(\beta(v^*)v)$. Then, 
\begin{eqnarray*}
ww^* & = & \alpha_s(\beta(v^*) vv^* \beta(v)) \\
& = & \alpha_s(\beta(v^*) \beta(vv^*) \beta(v)) \\
& = & \alpha_s \beta(v^*v) \neq 0.
\end{eqnarray*}
Hence, $w$ is a non-zero partial isometry in $q\widetilde{M}q$. Moreover, for every $x \in Q_0$,
\begin{eqnarray*}
w \alpha_{2s}(x) & = & \alpha_s(\beta(v^*) v \alpha_s(x)) \\
& = & \alpha_s(\beta(v^*) x v) \\
& = & \alpha_s(\beta(v^*x)v) \\
& = & \alpha_s(\beta(\alpha_s(x)v^*)v) \\
& = & \alpha_s\beta\alpha_s(x) \alpha_s(\beta(v^*)v) \\
& = & \beta(x) w \\
& = & xw.
\end{eqnarray*}

Since by induction, we can go till $s = 1$, we have found a non-zero partial isometry $v \in q\widetilde{M}q$ such that
\begin{equation}\label{inter}
xv = v\alpha_1(x), \forall x \in Q_0.
\end{equation}
Note that $v^*v \in \alpha_1(Q_0)' \cap qMq$. Moreover, since $\alpha_1 : q\widetilde{M}q \to q\widetilde{M}q$ is a $\ast$-automorphism, and $Q_0 \npreceq_{M} L(\R)$,  Theorem $2.4$ in \cite{houdayer4} gives
\begin{eqnarray*}
\alpha_1(Q_0)' \cap q\widetilde{M}q & = & \alpha_1\left( Q_0' \cap q\widetilde{M}q \right) \\
& \subset & \alpha_1(qMq).
\end{eqnarray*}
Hence $v^*v \in \alpha_1(qMq)$.

Since $Q_0 \npreceq_M L(\R)$, we know that there exists a sequence of unitaries $(u_k)$ in $Q_0$ such that $E_{L(\R)}(x^* u_k y) \to 0$ strongly, as $k \to \infty$, for any $x, y \in qM$. We need to go further and prove the following:

\begin{claim}\label{esperance}
$\forall a, b \in q\widetilde{M}q, \left\| E_{M_2}(a^* u_k b) \right\|_{2, \Tr} \to 0$, as $k \to \infty$.
\end{claim}

\begin{proof}[Proof of Claim $\ref{esperance}$]
Let $a, b \in (\widetilde{M})_1$ be either elements in $L(\R)$ or reduced words with letters alternating from $M_1 \ominus L(\R)$ and $M_2 \ominus L(\R)$. Write $b = y b'$ with
\begin{itemize}
\item $y = b$ if $b \in L(\R)$;
\item $y = 1$ if $b$ is a reduced word beginning with a letter from $M_2 \ominus L(\R)$;
\item $y =$ the first letter of $b$ otherwise. 
\end{itemize}
Note that either $b' = 1$ or $b'$ is a reduced word beginning with a letter from $M_2 \ominus L(\R)$.  Likewise write $a = a' x$ with
\begin{itemize}
\item $x = a$ if $x \in L(\R)$;
\item $x = 1$ if $a$ is a reduced word ending with a letter from $M_2 \ominus L(\R)$;
\item $x =$ the last letter of $a$ otherwise. 
\end{itemize}
Either $a' = 1$ or $a'$ is a reduced word ending with a letter from $M_2 \ominus L(\R)$. For any $z \in Q_0 \subset M_1$, $xzy - E_{L(\R)}(xzy) \in M_1 \ominus L(\R)$, so that
\begin{equation*}
E_{M_2} (a z b) = E_{M_2}(a' E_{L(\R)}(x z y) b').
\end{equation*}
Since $E_{L(\R)}(x u_k y) \to 0$ strongly, as $k \to \infty$, it follows that $E_{M_2}(a u_k b) \to 0$ strongly, as $k \to \infty$, as well. Thus, in the finite von Neumann algebra $q\widetilde{M}q$, we get $\left\| qE_{M_2}(a u_k b)q \right\|_{2, \Tr} \to 0$, as $k \to \infty$. 

Note that
\begin{equation*}
\mathcal{A}:= \mbox{span} \left\{ L(\R), (M_{i_1} \ominus L(\R)) \cdots (M_{i_n} \ominus L(\R)) : n \geq 1, i_1 \neq \cdots \neq i_n \right\}
\end{equation*}
is a unital $\ast$-strongly dense $\ast$-subalgebra of $\widetilde{M}$. What we have shown so far is that for any $a, b \in \mathcal{A}$, $\left\| qE_{M_2}(a u_k b)q \right\|_{2, \Tr} \to 0$, as $k \to \infty$. Let now $a, b \in (\widetilde{M})_1$. By Kaplansky density theorem, let $(a_i)$ and $(b_j)$ be sequences in $(\mathcal{A})_1$ such that $a_i \to a$ and $b_j \to b$ strongly. Recall that $(u_k)$ is a sequence in $Q_0 \subset q \widetilde{M} q$. We have
\begin{eqnarray*}
\left\| qE_{M_2}(a u_k b)q \right\|_{2, \Tr} & \leq & \left\| qE_{M_2}(a_i u_k b_j)q \right\|_{2, \Tr} + \left\| qE_{M_2}(a_i u_k (b - b_j))q \right\|_{2, \Tr} \\
& & + \left\| qE_{M_2}((a - a_i) u_k b_j)q \right\|_{2, \Tr} + \left\| qE_{M_2}((a - a_i) u_k (b - b_j))q \right\|_{2, \Tr} \\
& \leq & \left\| qE_{M_2}(a_i u_k b_j)q \right\|_{2, \Tr} + \left\| qa_i u_k (b - b_j)q \right\|_{2, \Tr} \\
& & + \left\| q(a - a_i) u_k b_j q \right\|_{2, \Tr} + \left\| q(a - a_i) u_k (b - b_j)q \right\|_{2, \Tr} \\
& \leq & \left\| qE_{M_2}(a_i u_k b_j)q \right\|_{2, \Tr} + 3 \left\|(b - b_j)q \right\|_{2, \Tr} + \left\| q(a - a_i)q\right\|_{2, \Tr}  \\
\end{eqnarray*}
Fix $\varepsilon > 0$. Since $a_i \to a$ and $b_j \to b$ strongly, let $i_0, j_0$ large enough such that 
\begin{equation*}
3 \left\|(b - b_{j_0})q \right\|_{2, \Tr} + \left\| q(a - a_{i_0})q\right\|_{2, \Tr} \leq \varepsilon/2.
\end{equation*}
Now let $k_0 \in \N$ such that for any $k \geq k_0$,  
\begin{equation*}
\left\| qE_{M_2}(a_{i_0} u_k b_{j_0})q \right\|_{2, \Tr} \leq \varepsilon/2.
\end{equation*}
We finally get $\left\| qE_{M_2}(a u_k b)q \right\|_{2, \Tr} \leq \varepsilon$, for any $k \geq k_0$, which finishes the proof of the claim. 
\end{proof}

We remind that for any $x \in Q_0$, $v^*xv = \alpha_1(x)v^*v$. Moreover, $v^*v \in \alpha_1(qMq) \subset qM_2q$. So, for any $x \in Q_0$, $v^*xv \in qM_2q$. Since $\alpha_1(u_k) \in \mathcal{U}(qM_2q)$, we get
\begin{equation*}
\left\| v^*v \right\|_{2, \Tr} = \left\| \alpha_1(u_k) v^*v \right\|_{2, \Tr} = \left\| E_{M_2}(\alpha_1(u_k) v^*v) \right\|_{2, \Tr} = \left\| E_{M_2}(v^* u_k v) \right\|_{2, \Tr} \to 0.
\end{equation*}
Thus $v = 0$, which is a contradiction. 
\end{proof}

\begin{proof}[Proof of Corollary $\ref{semisoliditybis}$]
Let $q \in L(\R)$ be a non-zero projection such that $\Tr(q) < \infty$. Denote by $N = qMq$ the corresponding ${\rm II_1}$ factor and by $\tau = \Tr(q \cdot q)/\Tr(q)$ the canonical trace on $N$. By contradiction, assume that $N$ is not semisolid. Then there exists $Q \subset N$ a non-amenable von Neumann subalgebra such that the relative commutant $Q' \cap N$ is of type ${\rm II_1}$. Write $z \in \mathcal{Z}(Q)$ for the maximal projection such that $Qz$ is amenable. Then $1 - z \neq 0$, the von Neumann algebra $Q(1 - z)$ has no amenable direct summand and $(Q' \cap N)(1 - z)$ is still of type ${\rm II_1}$. We may choose a projection $q_0 \in Q(1 - z)$ such that $\tau(q_0) = 1/n$. Since $N$ is a ${\rm II_1}$ factor, we may replace $Q$ by $\mathbf{M}_n(\C) \otimes q_0 Q q_0$, so that we may assume $Q \subset N$ has no amenable direct summand and $Q' \cap N$ is still of type ${\rm II_1}$.

If we apply Theorem \ref{semisolidity}, it follows that $Q' \cap N \preceq_M L(\R)$. We get a contradiction because $Q' \cap N$ is of type ${\rm II_1}$ and $L(\R)$ is of type ${\rm I}$.
\end{proof}

It follows from \cite{shlya2004} that for any type ${\rm III_1}$ factor $\mathcal{M}$, if the continuous core $M = \mathcal{M} \rtimes_\sigma \R$ is full, then Connes' invariant $\tau(\mathcal{M})$ is the usual topology on $\R$. Let now $\mathcal{M} = \Gamma(H_\R, U_t)''$ be a free Araki-Woods factor associated with $(U_t)$ an almost periodic representation. Denote by $S_U \subset \R^*_+$ the (countable) subgroup generated by the point spectrum of $(U_t)$. Then $\tau(\mathcal{M})$ is strictly weaker than the usual topology. More precisely, the completion of $\R$ w.r.t. the topology $\tau(\mathcal{M})$ is the compact group $\widehat{S_U}$ dual of $S_U$ (see \cite{connes74}). Therefore in this case, for any non-zero projection $q \in L(\R)$ such that $\Tr(q) < \infty$, the ${\rm II_1}$ factor $qMq$ is semisolid, by Theorem \ref{semisolidity}, and has property $\Gamma$ of Murray \& von Neumann by the above remark.

\subsection{Solidity of the continuous core under the assumption that $(U_t)$ is mixing}

We start this subsection with the following observations. The {\it solidity} of the continuous core $M$ forces the centralizers on $\mathcal{M}$ to be {\it amenable}. Indeed, fix $\psi$ any f.n. state on $\mathcal{M}$. Assume that the continuous core $M \simeq \mathcal{M} \rtimes_{\sigma^\psi} \R$ is solid. Choose
  a non-zero projection $q \in L(\R)$ such that $\Tr(q) < \infty$. Since $L(\R)q$ is diffuse in $q(\mathcal{M} \rtimes_{\sigma^\psi} \R)q$, its relative commutant must be amenable. In particular $\mathcal{M}^\psi \bar{\otimes} L(\R)q$ is amenable. Thus, $\mathcal{M}^\psi$ is amenable.

Note that if the orthogonal representation $(U_t)$ contains a $\frac{2\pi}{|\log \lambda|}$-periodic subrepresentation $(V_t^\lambda)$, $0 < \lambda < 1$, of the form
\begin{equation*}
V_t^\lambda = \begin{pmatrix}
\cos(t\log\lambda) & -\sin(t\log\lambda) \\
\sin(t\log\lambda) & \cos(t\log\lambda)
\end{pmatrix},
\end{equation*}
then the free Araki-Woods factor $\mathcal{M} = \Gamma(H_\R, U_t)''$ freely absorbs $L(\F_\infty)$ (see  \cite{shlya97}):
\begin{equation*}
(\mathcal{M}, \varphi_U) \ast (L(\F_\infty), \tau) \cong (\mathcal{M}, \varphi_U).
\end{equation*}
In particular, the centralizer of the free quasi-free state $\mathcal{M}^{\varphi_U}$ is non-amenable since it contains $L(\F_\infty)$. Therefore, whenever $(U_t)$ contains a periodic subrepresentation of the form $(V_t^\lambda)$ for some $0 < \lambda < 1$, the continuous core of $\Gamma(H_\R, U_t)''$ is semisolid by Theorem \ref{semisolidity} but can never be solid. However, when $(U_t)$ is assumed to be {\it mixing}, we get solidity of the continuous core. Indeed in that case, we can control the relative commutant $A' \cap M$ of diffuse subalgebras $A \subset L(\R) \subset M$, where $M$ is the continuous core of the  free Araki-Woods factor associated with $(U_t)$. Thus, the next theorem can be regarded as the analog of a result of Popa, namely Theorem 3.1 of \cite{popamal1} (see also Theorem D.4 in \cite{vaesbern}).

\begin{theo}\label{solidity}
Let $(U_t)$ be a mixing orthogonal representation of $\R$ on the real Hilbert space $H_\R$. Denote by $\mathcal{M} = \Gamma(H_\R, U_t)''$ the corresponding free Araki-Woods factor and by $M = \mathcal{M} \rtimes_\sigma \R$ its continuous core. Let $k \geq 1$ and let $q \in \mathbf{M}_k(\C) \otimes L(\R)$ be a non-zero projection such that $T := (\Tr_k \otimes \Tr)(q) < \infty$. Write $L(\R)^T := q(\mathbf{M}_k(\C) \otimes L(\R))q$ and $M^T := q(\mathbf{M}_k(\C) \otimes M)q$. Let $A \subset L(\R)^T$ be a diffuse von Neumann subalgebra. 

Then for any sub-bimodule $\vphantom{}_AH_{L(\R)^T}$ of $\vphantom{}_AL^2(M^T)_{L(\R)^T}$ such that $\dim(H_{L(\R)^T}) < \infty$, one has $H \subset L^2(L(\R)^T)$. In particular $A' \cap M^T \subset L(\R)^T$.
\end{theo}

\begin{cor}\label{soliditybis}
Let $\mathcal{M} = \Gamma(H_\R, U_t)''$ be a free Araki-Woods factor such that the orthogonal representation $(U_t)$ is mixing. Then the continuous core $M = \mathcal{M} \rtimes_\sigma \R$ is a solid ${\rm II_\infty}$ factor.
\end{cor}

\begin{proof}[Proof of Theorem $\ref{solidity}$]
As usual, denote by $(\lambda_t)$ the unitaries in $L(\R)$ that implement the modular action $\sigma$ on $\mathcal{M}$. Let $\Phi : L^\infty(\R) \to L(\R)$ be the Fourier Transform so that $\Phi(e^{it \cdot}) = \lambda_t$, for every $t \in \R$. Let $T > 0$ and denote by $q = \Phi(\chi_{[0, T]})$. Notice that $L^\infty(\R)\chi_{[0, T]} \cong L^\infty[0, T]$ and that
\begin{equation*}
\mbox{span} \left\{ \sum_{k \in F} c_k e^{i \frac{2\pi}{T} k \cdot}\chi_{[0, T]} :  F \subset \Z \mbox{ finite subset}, c_k \in \C, \forall k \in F \right\}
\end{equation*}
is a unital $\ast$-strongly dense $\ast$-subalgebra of $L^\infty(\R)\chi_{[0, T]}$. Thus, using the isomorphism $\Phi$, we get that 
\begin{equation*}
\mathcal{A} := \mbox{span} \left\{ \sum_{k \in F} c_k \lambda_{\frac{2\pi}{T}k} q : F \subset \Z \mbox{ finite subset}, c_k \in \C, \forall k \in F \right\}
\end{equation*}
is a unital $\ast$-strongly dense $\ast$-subalgebra of $L(\R)q$. Let $(u_n)$ be bounded sequence in $L(\R)q$ such that $u_n \to 0$ weakly, as $n \to \infty$, and $\left\| u_n \right\|_\infty \leq 1$, for every $n \in \N$. Using Kaplansky density theorem together with a standard diagonal process, choose a sequence $y_n \in \mathcal{A}$ such that $\left\| y_n \right\|_\infty \leq 1$, for every $n \in \N$, and $\left\| u_n - y_n \right\|_{2, \Tr} \to 0$, as $n \to \infty$. We will write $y_n = z_n q$ with
\begin{equation*}
z_n = \sum_{k \in F_n} c_{k, n} \lambda_{\frac{2 \pi}{T} k},
\end{equation*}
where $F_n \subset \Z$ is finite, $c_{k, n} \in \C$, for any $k \in F_n$ and any $n \in \N$. Using the $T$-periodicity, we have for any $n \in \N$,
\begin{eqnarray*}
\left\| z_n \right\|_\infty & = & \left\| \Phi^{-1}(z_n) \right\|_\infty \\
& = & \esssup_{x \in \R} \left| \sum_{k \in F_n} c_{k, n} e^{i \frac{2\pi}{T} k x} \right| \\
& = & \esssup_{x \in [0, T]} \left| \sum_{k \in F_n} c_{k, n} e^{i \frac{2\pi}{T} k x} \right| \\
& = & \left\| \Phi^{-1}(z_n) \chi_{[0, T]} \right\|_\infty \\
& = & \left\| y_n \right\|_\infty \leq 1.
\end{eqnarray*}
Thus, the sequence $(z_n)$ is uniformly bounded.

The {\bf first step} of the proof consists in proving the following:
\begin{equation*}
\left\| E_{L(\R)q}(a u_n b) \right\|_{2, \Tr} \to 0, \mbox{ as } n \to \infty, \forall a, b \in qMq \cap \ker\left(E_{L(\R)q}\right).
\end{equation*}
Equivalently, we need to show that
\begin{equation}\label{estimate}
\left\| q E_{L(\R)}(a u_n b) q \right\|_{2, \Tr} \to 0, \mbox{ as } n \to \infty, \forall a, b \in \ker\left(E_{L(\R)}\right).
\end{equation}

The {\bf first step} of the proof is now divided in three different claims that will lead to proving $(\ref{estimate})$. First note that 
\begin{equation*}
\mathcal{E} := \mbox{span}\left\{ \sum_{t \in F} x_t \lambda_{t} : F \subset \R \mbox{ finite subset}, x_t \in \mathcal{M} \mbox{ with } \varphi(x_t) = 0, \forall t \in F\right\}
\end{equation*}
is $\ast$-strongly dense in $\ker\left(E_{L(\R)}\right)$ by Kaplansky density theorem. We first prove the following:

\begin{claim}\label{estimate1}
If $\left\| q E_{L(\R)}(x u_n y) q \right\|_{2, \Tr} \to 0, \mbox{ as } n \to \infty$, $\forall x, y \in \mathcal{M}$ with $\varphi(x) = \varphi(y) = 0$, then $(\ref{estimate})$ is satisfied.
\end{claim}

\begin{proof}[Proof of Claim $\ref{estimate1}$]
Assume $\left\| q E_{L(\R)}(x u_n y) q \right\|_{2, \Tr} \to 0, \mbox{ as } n \to \infty$, $\forall x, y \in \mathcal{M}$ with $\varphi(x) = \varphi(y) = 0$. First take $a \in \mathcal{E}$ that we write $a = \sum_{s \in F} x_s \lambda_s$, with $F \subset \R$ finite subset, such that $x_s \in\mathcal{M}$, $\varphi(x_s) = 0$, for every $s \in F$. Then take $b \in \ker\left(E_{L(\R)}\right)$ and let $(b_j)_{j \in J}$ be a sequence in $\mathcal{E}$ such that $b - b_j \to 0$ $\ast$-strongly, as $j \to \infty$. Since $\left\| u_n \right\|_\infty \leq 1$, we get for any $n \in \N$ and any $j \in J$,
\begin{eqnarray*}
\left\| q E_{L(\R)}(a u_n b) q \right\|_{2, \Tr} & \leq & \left\| q E_{L(\R)}(a u_n b_j) q \right\|_{2, \Tr} + \left\| q E_{L(\R)}(a u_n (b - b_j)) q \right\|_{2, \Tr} \\
& \leq & \left\| q E_{L(\R)}(a u_n b_j) q \right\|_{2, \Tr} + \left\| a \right\|_\infty \left\| (b - b_j)q \right\|_{2, \Tr}
\end{eqnarray*} 
Fix $\varepsilon > 0$. Since $b - b_j \to 0$ $\ast$-strongly, as $j \to \infty$, fix $j_0 \in J$ such that $\left\| a \right\|_\infty \left\| (b - b_{j_0})q \right\|_{2, \Tr} \leq \varepsilon/2$. Write $b_{j_0} = \sum_{t \in F'} y_t \lambda_t$, with $F' \subset \R$ finite subset, such that $y_t \in \mathcal{M}$, $\varphi(y_t) = 0$, for every $t \in F'$. Therefore, for any $n \in \N$,
\begin{eqnarray*}
\left\| q E_{L(\R)}(a u_n b_{j_0}) q \right\|_{2, \Tr} & \leq & \sum_{(s, t) \in F \times F'} \left\| q E_{L(\R)}(x_s \lambda_s u_n y_t \lambda_t) q \right\|_{2, \Tr} \\
& = & \sum_{(s, t) \in F \times F'} \left\| \lambda_s q E_{L(\R)}(\sigma_{-s}(x_s) u_n y_t) q \lambda_t \right\|_{2, \Tr} \\
& = & \sum_{(s, t) \in F \times F'} \left\| q E_{L(\R)}(\sigma_{-s}(x_s) u_n y_t) q \right\|_{2, \Tr}.
\end{eqnarray*}
Since $\varphi(\sigma_{-s}(x_s)) = \varphi(y_t) = 0$, for any $(s, t) \in F \times F'$, using the assumption of the claim, there exists $n_0 \in \N$ large enough such that for any $n \geq n_0$, $\left\| q E_{L(\R)}(a u_n b_{j_0}) q \right\|_{2, \Tr} \leq \varepsilon/2$. Thus, for any $n \geq n_0$, $\left\| q E_{L(\R)}(a u_n b) q \right\|_{2, \Tr} \leq \varepsilon$. This proves that for any $a \in \mathcal{E}$ and any $b \in \ker\left(E_{L(\R)}\right)$, $\left\| q E_{L(\R)}(a u_n b) q \right\|_{2, \Tr} \to 0$, as $nÊ\to \infty$. If we do the same thing by approximating $a \in \ker\left(E_{L(\R)}\right)$ with elements in $\mathcal{E}$, using the fact that $u_n \in (L(\R)q)_1$, we finally get the claim.
\end{proof}

We now replace the sequence $(u_n)$ by $(z_n)$, use the mixing property of the modular action $\sigma$ and prove the following:

\begin{claim}\label{estimate2}
$\forall a, b \in (\mathcal M)_1$ with $\varphi(a) = \varphi(b) = 0$, $\left\| qE_{L(\R)}(a z_n b)q \right\|_{2, \Tr} \to 0$, as $n \to \infty$.
\end{claim}

\begin{proof}[Proof of Claim $\ref{estimate2}$]
Fix $a, b \in (\mathcal{M})_1$ such that $\varphi(a) = \varphi(b) = 0$. Fix $\varepsilon > 0$. For any $n \in \N$, we have
\begin{eqnarray*}
\left\| qE_{L(\R)}(a z_n b)q \right\|^2_{2, \Tr} & = &  \left\| \sum_{k \in F_n} c_{k, n} \varphi\left(a \sigma_{\frac{2\pi}{T} k}(b)\right) \lambda_{\frac{2\pi}{T} k} q \right\|^2_{2, \Tr} \\
& = & \Tr(q) \sum_{k \in F_n} |c_{k, n}|^2 \left|\varphi \left(a \sigma_{\frac{2\pi}{T} k}(b)\right)\right|^2.
\end{eqnarray*}
Moreover for any $n \in \N$,
\begin{equation*}
\Tr(q) \sum_{k \in F_n} |c_{k, n}|^2 = \left\| z_n q \right\|_{2, \Tr}^2 \leq \Tr(q) \left\| z_n q \right\|^2_\infty \leq T.
\end{equation*}
Since the modular group $\sigma$ is $\varphi$-mixing (because $(U_t)$ is assumed to be mixing), there exists a finite subset $K \subset \Z$ such that for any $k \in \Z \backslash K$, $\left|\varphi \left(a \sigma_{\frac{2\pi}{T} k}(b)\right)\right| \leq \varepsilon/\sqrt{2 T}$. Thus,
\begin{equation*}
\left\| qE_{L(\R)}(a z_n b)q \right\|_{2, \Tr} \leq \left\| \sum_{k \in K \cap F_n} c_{k, n} \lambda_{\frac{2\pi}{T} k} q \right\|_{2, \Tr} + \varepsilon/2.
\end{equation*}
Since $u_n - z_n q \to 0$ strongly and $u_n \to 0$ weakly, as $n \to \infty$, it follows that $z_n q \to 0$ weakly, as $n \to \infty$. In particular there exists $n_0$ large enough such that for any $n \geq n_0$, for any $k \in K \cap F_n$, $|c_{k, n}| \leq \varepsilon/(2 |K| \left\| q \right\|_{2, \Tr})$. Thus, for any $n \geq n_0$,
\begin{equation*}
\left\| qE_{L(\R)}(a z_n b)q \right\|_{2, \Tr} \leq \varepsilon/2 + \varepsilon/2 = \varepsilon.
\end{equation*}
This proves that $\left\| qE_{L(\R)}(a z_n b)q \right\|_{2, \Tr} \to 0$, as $n \to \infty$.
\end{proof}

The last claim consists in going back to the sequence $(u_n)$ and proving the following:

\begin{claim}\label{estimate3}
$\forall a, b \in (\mathcal M)_1$ with $\varphi(a) = \varphi(b) = 0$, $\left\| qE_{L(\R)}(a u_n b)q \right\|_{2, \Tr} \to 0$, as $n \to \infty$.
\end{claim}

\begin{proof}[Proof of Claim $\ref{estimate3}$]
Applying once more Kaplansky density theorem, we can find a sequence $(q_i)_{i \in I}$ in $L(\R)$ such that
\begin{itemize}
\item $q_i = \sum_{t \in F_i} d_t \lambda_t$, with $F_i \subset \R$ finite subset, $d_t \in \C$, for any $t \in F_i$ and for any $i \in I$;
\item $\|q_i\|_\infty \leq 1$, for any $i \in I$;
\item $q - q_i \to 0$ $\ast$-strongly, as $i \to \infty$.
\end{itemize}
Fix now $a, b \in (\mathcal M)_1$ such that $\varphi(a) = \varphi(b) = 0$. Using the fact that 
\begin{equation*}
\left\| a \right\|_\infty, \left\| b \right\|_\infty, \left\| q \right\|_\infty, \left\| z_n \right\|_\infty \leq 1, \forall n \in \N,
\end{equation*}
we get for any $n \in \N$ and any $i \in I$,
\begin{eqnarray*}
\left\| qE_{L(\R)}(a u_n b)q \right\|_{2, \Tr} & \leq & \left\| qE_{L(\R)}(a (u_n - z_n q) b)q \right\|_{2, \Tr} + \left\| qE_{L(\R)}(a z_n q b)q \right\|_{2, \Tr} \\
& \leq & \|u_n - z_n q\|_{2, \Tr} + \left\| qE_{L(\R)}(a z_n (q - q_i) b)q \right\|_{2, \Tr}  \\ & & + \left\| qE_{L(\R)}(a z_n q_i b)q \right\|_{2, \Tr} \\
& \leq & \|u_n - z_n q\|_{2, \Tr} + \left\| (q - q_i) b q \right\|_{2, \Tr}  \\ & & + \sum_{t \in F_i} |d_t| \left\| qE_{L(\R)}(a z_n \sigma_t(b)) \lambda_t q \right\|_{2, \Tr} \\
& \leq & \|u_n - z_n q\|_{2, \Tr} + \left\| (q - q_i) b q \right\|_{2, \Tr}  \\ & & + \sum_{t \in F_i} |d_t| \left\| qE_{L(\R)}(a z_n \sigma_t(b))q \right\|_{2, \Tr}.
\end{eqnarray*}
Since $q - q_i \to 0$ $\ast$-strongly, as $i \to \infty$, it follows that $\|(q - q_i)bq\|_{2, \Tr} \to 0$, as $i \to \infty$. Fix $\varepsilon > 0$. Then, take $i_0 \in I$ such that $\|(q - q_{i_0})bq\|_{2, \Tr} \leq \varepsilon/3$. Since $\|u_n - z_nq\|_{2, \Tr} \to 0$, as $n \to \infty$ and using Claim $\ref{estimate2}$, we may choose $n_0$ large enough such that for any $n \geq n_0$, 
\begin{eqnarray*}
\left\| u_n - z_n q \right\|_{2, \Tr} & \leq & \varepsilon/3 \\
\sum_{t \in F_{i_0}} |d_t| \left\| qE_{L(\R)}(a z_n \sigma_t(b))q \right\|_{2, \Tr} & \leq & \varepsilon/3.
\end{eqnarray*}
Consequently, for any $n \geq n_0$,  we get $\left\| qE_{L(\R)}(a u_n b)q \right\|_{2, \Tr} \leq \varepsilon$. Therefore, we have proven $\left\| qE_{L(\R)}(a u_n b)q \right\|_{2, \Tr} \to 0$, as $n \to \infty$.
\end{proof}

Thanks to Claims \ref{estimate1} and \ref{estimate3}, it is then clear that $(\ref{estimate})$ is satisfied. This finishes the {\bf first step} of the proof.

The {\bf last step} of the proof consists in using Theorem \ref{principle}. Let $k \geq 1$ and $q \in \mathbf{M}_k(\C) \otimes L(\R)$ be a non-zero projection such that $T := (\Tr_k \otimes \Tr)(q) < \infty$. Since $\mathbf{M}_k(\C) \otimes M$ is a ${\rm II_\infty}$ factor, there exists a unitary $u \in \mathcal{U}(\mathbf{M}_k(\C) \otimes M)$ such that
\begin{equation*}
q = 
u \begin{pmatrix}
q_0 & & 0 \\
& \ddots & \\
0 & & q_0
\end{pmatrix} u^*
\end{equation*}
where $q_0 = \Phi(\chi_{[0, T/k]}) \in L(\R)$.
Using the spatiality of $\Ad(u)$ on $\mathbf{M}_k(\C) \otimes M$, we may assume without loss of generality that 
\begin{equation*}
q = 
\begin{pmatrix}
q_0 & & 0 \\
& \ddots & \\
0 & & q_0
\end{pmatrix}
\end{equation*}
In particular, $q \in \mathbf{M}_k(\C) \otimes L(\R)q_0$. Define $M^T := q(\mathbf{M}_k(\C) \otimes M)q$ and $L(\R)^T := q(\mathbf{M}_k(\C) \otimes L(\R))q$. Let $A \subset L(\R)^T$ be a diffuse von Neumann subalgebra. Choose a sequence of unitaries $(u_n)$ in $A$ such that $u_n \to 0$ weakly, as $n \to \infty$. Thus, we can write $u_n = [u_n^{i,j}]_{i, j}$ where $u_n^{i, j} \in L(\R)q_0$ and $\left\| u_n^{i, j} \right\|_\infty \leq 1$, for any $n \in \N$ and any $i, j \in \{1, \dots, k\}$. Moreover, $u^{i, j}_n \to 0$ weakly, as $n \to \infty$, in $L(\R) q_0$, for any $i, j \in \{1, \dots, k\}$. Thus, using the {\bf first step} of the proof, it becomes clear that the inclusion $L(\R)^T \subset M^T$ is weakly mixing through $A$ in the sense of Definition \ref{weakmixing}. Thus, using Theorem \ref{principle}, it follows that for any $\vphantom{}_AH_{L(\R)^T}$ sub-bimodule of $\vphantom{}_AL^2(M^T)_{L(\R)^T}$ such that $\dim(H_{L(\R)^T}) < \infty$, one has $H \subset L^2(L(\R)^T)$. In particular $A' \cap M^T \subset L(\R)^T$.
\end{proof}

\begin{proof}[Proof of Corollary $\ref{soliditybis}$]
Let $q \in L(\R)$ be a non-zero projection such that $\Tr(q) < \infty$. Denote by $N = qMq$ the corresponding ${\rm II_1}$ factor. By contradiction assume that $N$ is not solid. Then there exists a non-amenable von Neumann subalgebra $Q \subset N$ such that the relative commutant $Q' \cap N$ is diffuse. Since $N$ is a ${\rm II_1}$ factor, using the same argument as in the proof of Corollary \ref{semisoliditybis}, we may assume that $Q$ has no amenable direct summand and $Q_0 = Q' \cap N$ is still diffuse.

Since $Q$ has no amenable direct summand, Theorem \ref{semisolidity} yields $Q_0 \preceq_M L(\R)$. Thus using Theorem \ref{intertwining}, we know that there exists a non-zero projection $p \in L(\R)$ such that $\Tr(p) < \infty$, and $Q_0 \preceq_{eMe} L(\R)p$ where $e = p \vee q$. Consequently, there exist $n \geq 1$, a (possibly non-unital) $\ast$-homomorphism $\psi : Q_0 \to \mathbf{M}_n(\C) \otimes L(\R)p$ and a non-zero partial isometry $v \in \mathbf{M}_{1, n}(\C) \otimes qMp$ such that 
\begin{equation*}
xv = v \psi(x), \forall x \in Q_0.
\end{equation*}
We moreover have 
\begin{equation*}
vv^* \in Q_0' \cap qMq \mbox{ and } v^*v \in \psi(Q_0)' \cap \psi(q)(\mathbf{M}_n(\C) \otimes pMp)\psi(q).
\end{equation*}
Write $Q_1 = Q_0' \cap qMq$ and notice that $Q \subset Q_1$. Since $\psi(Q_0)$ is diffuse and $v^*v \in \psi(Q_0)' \cap \psi(q)(\mathbf{M}_n(\C) \otimes pMp)\psi(q)$, Theorem \ref{solidity} yields $v^*v \in \psi(q)(\mathbf{M}_n(\C) \otimes L(\R)p)\psi(q)$, so that we may assume $v^*v = \psi(q)$. For any $y \in Q_1$, and any $x \in Q_0$,
\begin{eqnarray*}
v^*y v \psi(x) & = & v^* y x v \\
& = & v^* x y v \\
& = & \psi(x) v^* y v.  
\end{eqnarray*}
Thus, $v^* Q_1 v \subset \psi(Q_0)' \cap v^*v (\mathbf{M}_n(\C) \otimes pMp) v^*v$. Since $\psi(Q_0)$ is diffuse, Theorem \ref{solidity} yields $v^* Q_1 v  \subset v^*v(\mathbf{M}_n(\C) \otimes L(\R)p)v^*v$. Since $Q$ has no amenable direct summand and $Q \subset Q_1$ is a unital von Neumann subalgebra, it follows that $Q_1$ has no amenable direct summand either. Thus the von Neumann algebra $vv^* Q_1 vv^*$ is non-amenable. But $\Ad(v^*) : vv^* M vv^* \to v^*v(\mathbf{M}_n(\C) \otimes pMp)v^*v$ is a $\ast$-isomorphism and 
\begin{equation*}
\Ad(v^*)(vv^* Q_1 vv^*) \subset v^*v(\mathbf{M}_n(\C) \otimes L(\R)p)v^*v.
\end{equation*}
Since $v^*v(\mathbf{M}_n(\C) \otimes L(\R)p)v^*v$ is of type ${\rm I}$, hence amenable, we get a contradiction.
\end{proof}

Since the left regular representation $(\lambda_t)$ of $\R$ acting on $L^2_\R(\R, \mbox{Lebesgue})$ is mixing, the continuous core $M$ of $\Gamma(L^2_\R(\R, \mbox{Lebesgue}), \lambda_t)''$ is solid. We partially retrieve a previous result of Shlyakhtenko \cite{shlya98} where he proved in this case that $M \cong L(\F_\infty) \bar{\otimes} \mathbf{B}(\ell^2)$, which is solid by \cite{ozawa2003}. We will give in the next section an example of a non-amenable solid ${\rm II_1}$ factor with full fundamental group which is not isomorphic to any interpolated free group factor $L(\F_t)$, for $1 < t \leq \infty$.

Note that the mixing property of the representation $(U_t)$ is not a necessary condition for the solidity of the continuous core $M$. Indeed, take $U_t = \Id \oplus \lambda_t$ on $H_\R = \R \oplus L^2_\R(\R, \mbox{Lebesgue})$. Then $(U_t)$ is not mixing, but the continuous core $M$ of $\Gamma(H_\R, U_t)''$ is still isomorphic to $L(\F_\infty) \bar{\otimes} \mathbf{B}(\ell^2)$ \cite{shlya2003}.

\section{Examples of solid ${\rm II_1}$ factors}\label{solidfactor}

\subsection{Probability measures on the real line and unitary representations of $\R$} Write $\lambda$ for the Lebesgue measure on the real line $\R$. Let $\mu$ be a {\it symmetric} (positive) probability measure on $\R$, i.e. $\mu(X) = \mu(-X)$, for any Borel subset $X \subset \R$. Consider the following unitary representation $(U_t^\mu)$ of $\R$ on $L^2(\R, \mu)$ given by:
\begin{equation}\label{representation}
(U_t^\mu f)(x) = e^{itx} f(x), \forall f \in L^2(\R, \mu), \forall t, x \in \R. 
\end{equation}
Define the Hilbert subspace of $L^2(\R, \mu)$
\begin{equation}\label{hilbert}
K_\R^\mu := \left\{ f \in L^2(\R, \mu) : f(x) = \overline{f(-x)}, \forall x \in \R \right\}.
\end{equation} 
Since $\mu$ is assumed to be symmetric, the restriction of the inner product to $K_\R$ is real-valued. Indeed, for any $f, g \in K_\R^\mu$,
\begin{eqnarray*}
\langle f, g\rangle & = & \int_\R f(x) \overline{g(x)} \, d\mu(x) \\
& = & \int_\R f(-x) \overline{g(-x)} \, d\mu(-x) \\
& = & \int_\R \overline{f(x)} g(x) \, d\mu(x) \\
& = & \overline{\langle f, g\rangle}.
\end{eqnarray*}
Moreover the representation $(U_t^\mu)$ leaves $K_\R^\mu$ globally invariant. Thus, $(U_t^\mu)$ restricted to $K_\R^\mu$ becomes an orthogonal representation. Define the {\it Fourier Transform} of the probability measure $\mu$ by:
\begin{equation*}
\widetilde{\mu}(t) = \int_\R e^{itx} \, d\mu(x), \forall t \in \R.
\end{equation*}
We shall identify $\widehat{\R}$ with $\R$ in the usual way, such that
\begin{equation*}
\widehat{f}(t) = \int_\R e^{i t x}f(x) \,d\lambda(x), \forall t \in \R, \forall f \in L^1(\R, \lambda).
\end{equation*}

\begin{prop}\label{mixingrepresentation}
Let $\mu$ be a symmetric probability measure on $\R$. Then
\begin{equation*}
(U_t^\mu) \mbox{ is mixing } \Longleftrightarrow \widetilde{\mu}(t) \to 0, \mbox{ as } |t| \to \infty.
\end{equation*}
\end{prop}

\begin{proof}
We prove both directions.

$\Longrightarrow$ Assume $(U_t^\mu)$ is mixing. Let $f = \mathbf{1}_\R \in L^2(\R, \mu)$ be the constant function equal to $1$. Then
\begin{eqnarray*}
\widetilde{\mu}(t) & = & \int_\R e^{itx} \, d\mu(x) \\
& = & \langle U_t^\mu f, f\rangle \to 0, \mbox{ as } |t| \to \infty.
\end{eqnarray*}

$\Longleftarrow$ Assume $\widetilde{\mu}(t) \to 0$, as $|t| \to \infty$. Let $f, g \in L^2(\R, \mu)$. Then $h := f \overline{g} \in L^1(\R, \mu)$. Since the set $\left\{f \in \CC_0(\R) : \widehat{f} \in L^1(\R, \lambda) \right\}$ is dense in $L^1(\R, \mu)$, we may choose a sequence $(h_n)$ in C$_0(\R)$ such that $\left\| h - h_n \right\|_{L^1(\R, \mu)} \to 0$, as $n \to \infty$, and $\widehat{h}_n \in L^1(\R, \lambda)$, for any $n \in \N$. Define 
\begin{eqnarray*}
\widetilde{h}(t) & = & \int_\R e^{itx} h(x) \, d\mu(x), \forall t \in \R \\
\widetilde{h}_n(t) & = & \int_\R e^{itx} h_n(x) \, d\mu(x), \forall t \in \R, \forall n \in \N.
\end{eqnarray*}
Since $\left\| h - h_n \right\|_{L^1(\R, \mu)} \to 0$, as $n \to \infty$, it follows that $\left\| \widetilde{h} - \widetilde{h}_n \right\|_\infty \to 0$, as $n \to \infty$. Since $\widehat{h}_n \in L^1(\R, \lambda)$, we know that
\begin{equation*}
h_n(x) = C \int_\R e^{-ixu} \widehat{h}_n(u) \, d\lambda(u), \forall x \in \R,
\end{equation*}
where $C$ is a universal constant that only depends on the normalization of the Lebesgue measure $\lambda$ on $\R$. Therefore, for any $t \in \R$ and any $n \in \N$,
\begin{eqnarray*}
\widetilde{h}_n(t) & = & \int_{x \in \R} e^{itx} h_n(x) \, d\mu(x) \\
& = & C \int_{x \in \R} \left( \int_{u \in \R} e^{i(t - u)x} \widehat{h}_n(u) \, d\lambda(u) \right) d\mu(x) \\
& = & C \int_{u \in \R} \widehat{h}_n(u) \left( \int_{x \in \R} e^{i(t - u)x} \, d\mu(x)\right) d\lambda(u) \\
& = & C \int_{u \in \R} \widehat{h}_n(u) \widetilde{\mu}(t - u) \, d\lambda(u) \\
& = & C \left(\widehat{h}_n \ast \widetilde{\mu}\right) (t),
\end{eqnarray*}
where $\ast$ is the convolution product. Since $\widetilde{\mu} \in \CC_0(\R)$ and $\widehat{h}_n \in L^1(\R, \lambda)$, it is easy to check that $\widehat{h}_n \ast \widetilde{\mu} \in \CC_0(\R)$. Consequently, $\widetilde{h}_n \in \CC_0(\R)$ and since $\left\| \widetilde{h} - \widetilde{h}_n \right\|_\infty \to 0$, as $n \to \infty$, it follows that $\widetilde{h} \in \CC_0(\R)$. But for any $t \in \R$,
\begin{eqnarray*}
\langle U_t^\mu f, g\rangle & = & \int_\R e^{itx} f(x) \overline{g(x)} \, d\mu(x) \\
& = & \widetilde{h}(t).
\end{eqnarray*}
Thus, the unitary representation $(U^\mu_t)$ is mixing.
\end{proof}

For a measure $\nu$ on $\R$, define the {\it measure class} of $\nu$ by:
\begin{equation*}
\mathcal{C}_\nu := \left\{ \nu' : \nu' \mbox{ is absolutely continuous w.r.t. } \nu \right\}.
\end{equation*}

\begin{df}
Let $(V_t)$ be a unitary representation of $\R$ on a separable Hilbert space $H$. Denote by $B$ the infinitesimal generator of $(V_t)$, i.e. $B$ is the positive, self-adjoint (possibly) unbounded operator on $H$ such that $V_t = B^{it}$, for every $t \in \R$. We define the {\it spectral measure} of the representation $(V_t)$ as the spectral measure of the operator $B$ and denote it by $\mathcal{C}_{V}$. 
\end{df}

The {\it measure class} $\mathcal{C}_V$ can also be defined as the smallest collection of all the  measures $\nu$ on $\R$ such that:
\begin{enumerate}
\item If $\nu \in \mathcal{C}_V$ and $\nu'$ is absolutely continuous w.r.t. $\nu$, then $\nu' \in \mathcal{C}_V$;
\item For any unit vector $\eta \in H$, the probability measure associated with the positive definite function $t \mapsto \langle V_t \eta, \eta\rangle$ belongs to $\mathcal{C}_V$.
\end{enumerate}
Since $H$ is separable, there exists a measure $\nu$ that generates $\mathcal{C}_V$, i.e. $\mathcal{C}_V$ is the smallest collection of measures on $\R$ satisfying $(1)$ and containing $\nu$. We will refer to this particular measure $\nu$ as the {\textquotedblleft spectral measure\textquotedblright} of the representation $(V_t)$ and simply denote it by $\nu$. 

Let $\mu$ be a symmetric probability measure on $\R$ and consider the unitary representation $(U_t^\mu)$ on $L^2(\R, \mu)$ as defined in $(\ref{representation})$. Then for any unit vector $f \in L^2(\R, \mu)$, 
\begin{equation*}
\langle U_t^\mu f, f\rangle = \int_\R e^{itx} \left| f(x) \right|^2 \, d\mu(x), \forall t \in \R.
\end{equation*}
Since the probability measure $|f(x)|^2 \, d\mu(x)$ is absolutely continuous w.r.t. $d\mu(x)$, it is clear that the spectral measure of $(U_t^\mu)$ is $\mu$. More generally, we have the following:

\begin{prop}\label{spectralmeasure}
Let $\mu$ be a symmetric probability measure on $\R$. Consider the unitary representation $(U^\mu_t)$ defined on $L^2(\R, \mu)$ by $(\ref{representation})$. Then for any $n \geq 1$, the spectral measure of the $n$-fold tensor product $(U^\mu_t)^{\otimes n}$ is the $n$-fold convolution product 
\begin{equation*}
\mu^{\ast  n} = \underbrace{\mu \ast \cdots \ast \mu}_{n \, \fois}.
\end{equation*}
\end{prop}

\subsection{Examples of solid ${\rm II_1}$ factors}

Erd\"os showed in \cite{erdos} that the symmetric probability measure $\mu_\theta$, with $\theta = 5/2$, obtained as the weak limit of
\begin{equation*}
\left( \frac12 \delta_{-\theta^{-1}} + \frac12 \delta_{\theta^{-1}} \right) \ast \cdots \ast \left( \frac12 \delta_{-\theta^{-n}} + \frac12 \delta_{\theta^{-n}} \right) 
\end{equation*}
has a Fourier Transform
\begin{equation*}
\widetilde{\mu}_\theta(t) = \prod_{n \geq 1} \cos\left(\frac{t}{\theta^n}\right)
\end{equation*}
which vanishes at infinity, i.e. $\widetilde{\mu}(t) \to 0$, as $|t| \to \infty$, and $\mu_\theta$ is singular w.r.t. the Lebesgue measure $\lambda$.

\begin{exam}\label{singularmeasure}
Modifying the measure $\mu_\theta$, Antoniou \& Shkarin (see Theorem $2.5, {\rm v}$ in \cite{antoniou}) constructed an example of a symmetric probability $\mu$ on $\R$ such that:
\begin{enumerate}
\item The Fourier Transform of $\mu$ vanishes at infinity, i.e. $\widetilde{\mu}(t) \to 0$, as $|t| \to \infty$.
\item For any $n \geq 1$, the $n$-fold convolution product $\mu^{\ast n}$ is singular w.r.t. the Lebesgue measure $\lambda$.
\end{enumerate}
\end{exam}

Let $\mu$ be a symmetric probability measure on $\R$ as in Example \ref{singularmeasure}. Proposition \ref{mixingrepresentation} and Proposition \ref{spectralmeasure} yields that the unitary representation $(U_t^\mu)$ defined on $L^2(\R, \mu)$ by $(\ref{representation})$ satisfies:
\begin{enumerate}
\item $(U_t^\mu)$ is mixing.
\item The spectral measure of $\bigoplus_{n \geq 1} (U_t^\mu)^{\otimes n}$ is singular w.r.t. the Lebesgue measure $\lambda$.
\end{enumerate} 

Let now $\mathcal{M} = \Gamma(H_\R, U_t)''$ and let $M = \mathcal{M} \rtimes_\sigma \R$ be the continuous core. Let $q \in L(\R)$ be a non-zero projection such that $\Tr(q) < \infty$. Denote by $N = qMq$ the corresponding ${\rm II_1}$ factor. Using free probability techniques such as the {\it free entropy}, Shlyakhtenko (see Theorem $9.12$ in \cite{shlya99}) showed that if the spectral measure of the unitary representation $\bigoplus_{n \geq 1} U_t^{\otimes n}$ is singular w.r.t. the Lebesgue measure $\lambda$, then for any finite set of generators $X_1, \dots, X_n$ of $N$, the free entropy dimension satisfies
\begin{equation*}
\delta_0(X_1, \dots, X_n) \leq 1.
\end{equation*}
In particular, $N$ is not isomorphic to any interpolated free group factor $L(\F_t)$, for $1 < t \leq \infty$. Combining these two results together with Corollary \ref{solidity}, we obtain the following:

\begin{theo}\label{examplesolidfactor}
Let $\mu$ be a symmetric probability measure on $\R$ as in Example $\ref{singularmeasure}$. Let $\mathcal{M} = \Gamma(K_\R^\mu, U^\mu_t)''$ be the free Araki-Woods factor associated with the orthogonal representation $(U_t^\mu)$ acting on the real Hilbert space $K_\R^\mu$, as defined in $(\ref{representation}-\ref{hilbert})$. Let $M = \mathcal{M} \rtimes_\sigma \R$ be the continuous core. Fix a non-zero projection $q \in L(\R)$ such that $\Tr(q) < \infty$, and denote by $N = qMq$ the corresponding ${\rm II_1}$ factor. Then
\begin{enumerate}
\item $N$ is non-amenable and solid.
\item $N$ has full fundamental group, i.e. $\mathcal{F}(N) = \R^*_+$.
\item $N$ is not isomorphic to any interpolated free group factor $L(\F_t)$, for $1 < t \leq \infty$.
\end{enumerate}
\end{theo}

We believe that all the free Araki-Woods factors $\mathcal{M} = \Gamma(H_\R, U_t)''$ have the {\it complete metric approximation property} (c.m.a.p.), i.e. there exists a sequence $\Phi_n : \mathcal{M} \to \mathcal{M}$ of finite rank, completely bounded maps such that $\Phi_n \to \Id$ ultraweakly pointwise, as $n \to \infty$, and $\limsup_{n \to \infty} \left\| \Phi_n \right\|_{\cb} \leq 1$. If $\mathcal{M} = \Gamma(H_\R, U_t)''$ had the c.m.a.p. then by \cite{anan95}, the continuous core $M = \mathcal{M} \rtimes_\sigma \R$ would have the c.m.a.p., as well as the ${\rm II_1}$ factor $qMq$, for $q \in M$ non-zero finite projection. On the other hand, the wreath product ${\rm II_1}$ factors $L(\Z \wr \F_n)$ do not have the c.m.a.p., for any $2 \leq n \leq \infty$, by \cite{ozawapopa}. Thus, we conjecture that the solid ${\rm II_1}$ factors constructed in Theorem \ref{examplesolidfactor} are not isomorphic to $L(\Z \wr \F_n)$, for any $2 \leq n \leq \infty$.

\bibliographystyle{plain}

\end{document}